\documentclass[11pt]{amsart}

\usepackage[all]{xy}
\usepackage{mathrsfs}
\usepackage{amsmath,amsthm, amssymb, amsgen,amsxtra,amsfonts, amsbsy} 
\usepackage[all]{xy}
\usepackage{verbatim}
\usepackage{enumerate}

\usepackage{graphicx}
\usepackage{caption,rotating}
\usepackage{color}
\usepackage{subcaption}
\usepackage{hyperref}
\usepackage{times}
\newcommand{\mathds}[1]{{\mathbb #1}}

\begin{document}
%
%   D e f i n i t i o n s
%
%
\theoremstyle{definition}
\newtheorem{Definition}{Definition}[section]
\newtheorem*{Definitionx}{Definition}
\newtheorem{Convention}{Convention}[section]
\newtheorem{Notation}[Definition]{Notation}
\newtheorem{Construction}{Construction}[section]
\newtheorem{Example}[Definition]{Example}
\newtheorem{Examples}[Definition]{Examples}
\newtheorem{Remark}[Definition]{Remark}
\newtheorem*{Remarkx}{Remark}
\newtheorem{Remarks}[Definition]{Remarks}
\newtheorem{Caution}[Definition]{Caution}
\newtheorem{Conjecture}[Definition]{Conjecture}
\newtheorem*{Conjecturex}{Conjecture}
\newtheorem{Question}[Definition]{Question}
\newtheorem*{Questionx}{Question}
\newtheorem*{Acknowledgements}{Acknowledgements}
\newtheorem*{Organization}{Organization}
\newtheorem*{Disclaimer}{Disclaimer}
\theoremstyle{plain}
\newtheorem{Theorem}[Definition]{Theorem}
\newtheorem*{Theoremx}{Theorem}
\newtheorem{Proposition}[Definition]{Proposition}
\newtheorem*{Propositionx}{Proposition}
\newtheorem{Lemma}[Definition]{Lemma}
\newtheorem{Corollary}[Definition]{Corollary}
\newtheorem*{Corollaryx}{Corollary}
\newtheorem{Fact}[Definition]{Fact}
\newtheorem{Facts}[Definition]{Facts}
\newtheoremstyle{voiditstyle}{3pt}{3pt}{\itshape}{\parindent}%
{\bfseries}{.}{ }{\thmnote{#3}}%
\theoremstyle{voiditstyle}
\newtheorem*{VoidItalic}{}
\newtheoremstyle{voidromstyle}{3pt}{3pt}{\rm}{\parindent}%
{\bfseries}{.}{ }{\thmnote{#3}}%
\theoremstyle{voidromstyle}
\newtheorem*{VoidRoman}{}

\newcommand{\prf}{\par\noindent{\sc Proof.}\quad}
\newcommand{\blowup}{\rule[-3mm]{0mm}{0mm}}
\newcommand{\cal}{\mathcal}
\newcommand{\Aff}{{\mathds{A}}}
\newcommand{\BB}{{\mathds{B}}}
\newcommand{\CC}{{\mathds{C}}}
\newcommand{\EE}{{\mathds{E}}}
\newcommand{\FF}{{\mathds{F}}}
\newcommand{\GG}{{\mathds{G}}}
\newcommand{\HH}{{\mathds{H}}}
\newcommand{\NN}{{\mathds{N}}}
\newcommand{\ZZ}{{\mathds{Z}}}
\newcommand{\PP}{{\mathds{P}}}
\newcommand{\QQ}{{\mathds{Q}}}
\newcommand{\RR}{{\mathds{R}}}

\newcommand{\p}         {\phantom{{}^*}}
\newcommand{\Liea}{{\mathfrak a}}
\newcommand{\Lieb}{{\mathfrak b}}
\newcommand{\Lieg}{{\mathfrak g}}
\newcommand{\Liem}{{\mathfrak m}}
\newcommand{\ideala}{{\mathfrak a}}
\newcommand{\idealb}{{\mathfrak b}}
\newcommand{\idealg}{{\mathfrak g}}
\newcommand{\idealm}{{\mathfrak m}}
\newcommand{\idealp}{{\mathfrak p}}
\newcommand{\idealq}{{\mathfrak q}}
\newcommand{\idealI}{{\cal I}}
\newcommand{\lin}{\sim}
\newcommand{\num}{\equiv}
\newcommand{\dual}{\ast}
\newcommand{\iso}{\cong}
\newcommand{\homeo}{\approx}
\newcommand{\mm}{{\mathfrak m}}
\newcommand{\pp}{{\mathfrak p}}
\newcommand{\qq}{{\mathfrak q}}
\newcommand{\rr}{{\mathfrak r}}
\newcommand{\pP}{{\mathfrak P}}
\newcommand{\qQ}{{\mathfrak Q}}
\newcommand{\rR}{{\mathfrak R}}
%
%  evtl. auch \"uber \mathbb oder \Bbb
%
\newcommand{\OO}{{\cal O}}
\newcommand{\numero}{{n$^{\rm o}\:$}}
\newcommand{\mf}[1]{\mathfrak{#1}}
\newcommand{\mc}[1]{\mathcal{#1}}
\newcommand{\into}{{\hookrightarrow}}
\newcommand{\onto}{{\twoheadrightarrow}}
\newcommand{\Spec}{{\rm Spec}\:}
\newcommand{\BigSpec}{{\rm\bf Spec}\:}
\newcommand{\Spf}{{\rm Spf}\:}
\newcommand{\Proj}{{\rm Proj}\:}
\newcommand{\Pic}{{\rm Pic }}
\newcommand{\Num}{{\rm Num }}
\newcommand{\MW}{{\rm MW }}
\newcommand{\Br}{{\rm Br}}
\newcommand{\NS}{{\rm NS}}
\newcommand{\CH}{{\rm CH}}
\newcommand{\Sym}{{\mathfrak S}}
\newcommand{\Aut}{{\rm Aut}}
\newcommand{\Autp}{{\rm Aut}^p}
\newcommand{\Hom}{{\rm Hom}}
\newcommand{\Ext}{{\rm Ext}}
\newcommand{\ord}{{\rm ord}}
\newcommand{\coker}{{\rm coker}\,}
\newcommand{\Char}{{\rm char}}
\newcommand{\Ker}{{\rm ker}}
 \newcommand{\im}{{\rm im}}
\newcommand{\divisor}{{\rm div}}
\newcommand{\Def}{{\rm Def}}
\newcommand{\piet}{{\pi_1^{\rm \acute{e}t}}}
\newcommand{\Het}[1]{{H_{\rm \acute{e}t}^{{#1}}}}
\newcommand{\Hfl}[1]{{H_{\rm fl}^{{#1}}}}
\newcommand{\Hcris}[1]{{H_{\rm cris}^{{#1}}}}
\newcommand{\HdR}[1]{{H_{\rm dR}^{{#1}}}}
\newcommand{\hdR}[1]{{h_{\rm dR}^{{#1}}}}
\newcommand{\loc}{{\rm loc}}
\newcommand{\et}{{\rm \acute{e}t}}
\newcommand{\defin}[1]{{\bf #1}}
\newcommand{\oX}{\cal{X}}
\newcommand{\oA}{\cal{A}}
\newcommand{\oY}{\cal{Y}}
\newcommand{\blue}{\textcolor{blue}}
\newcommand{\red}{\textcolor{red}}

\author{Gebhard Martin}
\address{\hfill \newline
Mathematisches Institut  \newline
Universit\"at Bonn \newline
Endenicher Allee 60 \newline
53115 Bonn \newline
Germany}
\email{gmartin@math.uni-bonn.de}

\title{Infinitesimal automorphisms of algebraic varieties\\ \hspace{-3mm}and vector fields on elliptic surfaces}

\date{\today}

\begin{abstract}
We give several results concerning the connected component $\Aut_X^0$ of the automorphism scheme of a proper variety $X$ over a field, such as its behaviour with respect to birational modifications, normalization, restrictions to closed subschemes and deformations. Then, we apply our results to study the automorphism scheme of not necessarily Jacobian elliptic surfaces $f: X \to C$ over algebraically closed fields, generalizing work of Rudakov and Shafarevich, while giving counterexamples to some of their statements. We bound the dimension $h^0(X,T_X)$ of the space of global vector fields on an elliptic surface $X$ if the generic fiber of $f$ is ordinary or if $f$ admits no multiple fibers, and show that, without these assumptions, the number $h^0(X,T_X)$ can be arbitrarily large for any base curve $C$ and any field of positive characteristic. If $f$ is not isotrivial, we prove that $\Aut_X^0 \cong \mu_{p^n}$ and give a bound on $n$ in terms of the genus of $C$ and the multiplicity of multiple fibers of $f$. As a corollary, we re-prove the non-existence of global vector fields on K3 surfaces and calculate the connected component of the automorphism scheme of a generic supersingular Enriques surface in characteristic $2$. Finally, we present additional results on horizontal and vertical group scheme actions on elliptic surfaces which can be applied to determine $\Aut_X^0$ explicitly in many concrete cases.
\end{abstract}

\maketitle

\vspace{-5mm}
{\bf Contents}
\begin{itemize}
\item[\S \ref{introduction}]
Introduction \hfill \pageref{introduction}

\item[\S \ref{preliminaries}] 
Generalities \hfill \pageref{preliminaries}

\item[\S \ref{vertical}]
Vertical component of $\Aut_X^0$ \hfill \pageref{vertical}

\item[\S \ref{horizontal}]  
Horizontal component of $\Aut_X^0$  \hfill \pageref{horizontal}

\item[\S \ref{examples}]
Examples \hfill \pageref{examples}

\item[\S \ref{proofs}]
Proofs of the main theorems \hfill \pageref{proofs}

\item[\S \ref{appendix}] 
Appendix: Some quadratic twists \hfill \pageref{appendix}
\end{itemize}

\section{Introduction} \label{introduction}
Let $X$ be a scheme which is proper over a field $k$. In  \cite{MatsumuraOort} Matsumura and Oort proved that the automorphism functor $\Aut_X$ of $X$ over $k$ is representable by a group scheme that is locally of finite type over $k$. Its connected component of the identity $\Aut_X^0$ together with its tangent space at the identity  $H^0(X,T_X)$, consisting of global vector fields, play a central r\^ole in the deformation theory of $X$. Indeed, if $h^0(X,T_X) = 0$, then the deformation funtor $\Def_X$ of $X$ is prorepresentable and, conversely, if $\Aut_X^0$ is not smooth, then $\Def_X$ can never be prorepresentable. Similarly, if $X$ is a proper variety with $h^0(X,T_X) =  0$ and if a moduli stack $\cal{M}$ parametrizing objects of the same type as $X$ exists, this stack is Deligne--Mumford at the point corresponding to $X$, since the stabilizer of $\cal{M}$ at $X$ is reduced. This leads to the following geometric question.

\begin{Questionx}[{\bf A}]\label{QuestionA}
Let $X$ be a proper scheme over $k$. What is the dimension of $H^0(X,T_X)$?
\end{Questionx}

On the other hand, it was observed in \cite{AbramovichOlssonVistoli} and \cite{Alper} that Artin stacks with finite linearly reductive stabilizers behave better than general Deligne--Mumford stacks in many ways; for example, they are \'etale locally quotient stacks by finite and linearly reductive group schemes (see \cite[Theorem 3.2]{AbramovichOlssonVistoli}). Since linearly reductive group schemes may very well be connected, a first step towards checking whether $\Aut_X$ is linearly reductive is to check it for $\Aut_X^0$.

\begin{Questionx}[{\bf B}]\label{QuestionB}
When is $\Aut_X^0$ linearly reductive?
\end{Questionx}

Finally, to better understand how close $\cal{M}$ is to being a scheme \'etale locally at $X$, we can ask for the size of $\Aut_X^0$.

\begin{Questionx}[{\bf C}]\label{QuestionC}
If $\Aut_X^0$ is finite, what is its length?
\end{Questionx}

If $X$ is a smooth projective curve, answers to all three of the above questions are known, since $\Aut_X^0$ is always smooth in this case and the automorphism groups of $X$ are well-known. More precisely, we have $\Aut_{\PP^1} \cong {\rm PGL}_2$ and $\Aut_E^0 \cong E$ for an elliptic curve $E$ and if $X$ has higher genus, then $\Aut_X^0$ is trivial. However, already for singular curves or smooth projective surfaces, the three Questions (A), (B), and (C) are wide open. In recent years, however, some structural results in the case of surfaces of general type have been obtained by Tziolas in \cite{Tziolas1} and \cite{Tziolas2}.

The purpose of this paper is to give answers to Questions (A), (B) and (C) for elliptic surfaces over an algebraically closed field $k$ of arbitrary characteristic. Before we start explaining our setup, let us remark that Question (A) for elliptic surfaces without multiple fibers has been studied by Rudakov and Shafarevich \cite{RudakovShafarevich} and the first proof of the non-existence of global vector fields on K3 surfaces is a corollary of their work. Unfortunately, it turns out that Lemma $3$, Lemma $4$, and, as a result, also Theorem $6$ in \cite{RudakovShafarevich} are false as they are stated there. Some of these issues were also addressed in \cite{RudaShafaVector}, but the classification of counterexamples stated there is incomplete. In Section \ref{examples}, we give counterexamples to these claims and complete their classification of counterexamples begun in \cite{RudaShafaVector}. As a special case of our analysis, we will recover a modified version of \cite[Theorem 6]{RudakovShafarevich} in Theorem (D), which gives a characterization of elliptic surfaces with vector fields and without multiple fibers.

Let $k$ be an algebraically closed field of characteristic $\Char(k) = p \geq 0$. Let $f: X \to C$ be an elliptic surface, that is, $X$ is a smooth projective surface, $C$ is a smooth projective curve and $f$ is a proper morphism such that $f_*\OO_X = \OO_C$, almost all fibers of $f$ are smooth curves of genus one and $f$ is relatively minimal, that is, there are no $(-1)$-curves in the fibers of $f$. By Blanchard's Lemma (see \cite[Theorem 7.2.1]{Brion}), there is a natural morphism of group schemes $f_*: \Aut_X^0 \to \Aut_C^0$. We say that $\Ker(f_*)$ is the \emph{vertical component} and $\im(f_*)$ is the \emph{horizontal component} of $\Aut_X^0$.

In characteristic $0$, the structure of elliptic surfaces with non-trivial $\Aut_X^0$ is simple and well-known: If $\Aut_X^0$ is non-trivial, then either $X$ is ruled or, after a finite base change $C' \to C$, it becomes isomorphic to the trivial elliptic surface $F \times C'$ where $F$ is a general fiber of $f$. We leave it to the reader to check that the same conclusion follows from our results in arbitrary characteristic under the stronger assumption $\dim \Aut_X^0 > 0$ (see also the very recent preprint by Fong \cite{Fong}).

Recall that an elliptic surface is called \emph{isotrivial} if all smooth fibers of $f$ are isomorphic, or equivalently, if the $j$-map of $f$ is constant. 
In the following, in Theorems (A), (B), and (C), we will give a summary of our answers to Questions (A), (B), and (C). We refer the reader to Section \ref{vertical} and Section \ref{horizontal} for more refined and more general results on the structure of $\Ker(f_*)$ and $\im(f_*)$, such as possible fiber types, further information on multiplicities of fibers, as well as geometric restrictions on multisections for elliptic surfaces with non-trivial $\Aut_X^0$. The proofs of the following results can be found in Section \ref{proofs}, where we combine our results on vertical and horizontal components of $\Aut_X^0$.

\begin{Theoremx}[{\bf A}]
Let $f: X \to C$ be an elliptic surface. Then, the following hold:
\begin{enumerate}[(i)]
\item If $f$ is not isotrivial, then $h^0(X,T_X) \leq 1$.
\item If the generic fiber of $f$ is ordinary or $f$ admits no multiple fibers, then $h^0(X,T_X) \leq 4$. If additionally $h^0(X,T_X) \geq 2$, then one of the following holds:
\begin{enumerate}[(1)]
\item $X$ is ruled over an elliptic curve.
\item $X$ is an Abelian surface isogeneous to a product of elliptic curves.
\item $X$ is a bielliptic surface with $\omega_X \cong \OO_X$. These surfaces exist if and only if $p \in \{2,3\}$.
\end{enumerate}
\item For every field $K$ of characteristic $\Char(K) > 0$, for every smooth projective curve $C$ over $K$ and for every $n \geq 0$, there is an elliptic surface $f: X \to C$ with $h^0(X,T_X) \geq n$.
\end{enumerate}
\end{Theoremx}

In particular, the elliptic surfaces appearing in Theorem (A) (iii) are isotrivial with supersingular generic fiber and they admit multiple fibers. The relevant examples can be found in Example \ref{importantexample}. In the non-isotrivial cases, the following theorems give a description of $\Aut_X^0$ as well as a bound on its length that depends on the number $h_p$, which is defined in the discussion before Proposition \ref{Igusainequality} and which coincides with the number of supersingular $j$-invariants over $k$ if $p \neq 2,3$.

\begin{Theoremx}[{\bf B}]
Let $f: X \to C$ be a non-isotrivial elliptic surface.  Then, $\Aut_X^0 \cong \mu_{p^n}$ for some $n \geq 0$. In particular, $\Aut_X^0$ is linearly reductive.
\end{Theoremx}

\begin{Theoremx}[{\bf C}]
Let $f: X \to C$ be a non-isotrivial elliptic surface with $\Aut_X^0 \cong \mu_{p^n}$. Then, 
$$
\frac{1}{48}(p-1)(p^{2n-1} - 12 p^{n-1} + 1) + 1 - \frac{h_p}{2} \leq g(C).
$$
If, additionally,
\begin{enumerate}[(a)]
\item $p^n \geq 4$, or
\item $C \not \cong \PP^1$, or
\item $p^n = 3$ and the singular fibers of $f$ are not of type $({\rm II},{\rm I}_{3^{2k}})$ or $({\rm II},{\rm I}^*_{3^{2k-1}})$ for any $k \geq 1$, or
\item $p^n = 2$ and the singular fibers of $f$ are not of type $({\rm II},{\rm I}_{2^{2k+1}})$ or $({\rm III},{\rm I}_{2^{2k+1}})$ for any $k \geq 1$,
\end{enumerate}
then $\Ker(f_*) \cong \Aut_X^0$ and all additive or supersingular fibers of $f$ are multiple fibers with multiplicity divisible by $p^n$.
\end{Theoremx}

The exceptions in Theorem (C) (c) and (d) occur for every $k \geq 1$ (see Example \ref{caseb3} and Example \ref{caseb2}). Our Examples \ref{caseb3} and Example \ref{caseb2} are elliptic surfaces with a section and therefore they are counterexamples to \cite[Theorem 6]{RudakovShafarevich} in characteristic $2$ and $3$. Three of these four families of counterexamples were already exhibited in \cite{RudaShafaVector}.
In the other characteristics, our analysis recovers \cite[Theorem 6]{RudakovShafarevich}. More precisely, we obtain the following theorem.

\begin{Theoremx}[{\bf D}]
Let $f: X \to C$ be an elliptic surface without multiple fibers and such that $h^0(X,T_X) \neq 0$. Then, one of the following holds:
\begin{enumerate}[(i)]
\item $f$ is isotrivial and $c_2(X) = 0$.
\item $f$ is Jacobian and isotrivial with two singular fibers, $X$ is rational, and $C \cong \PP^1$.
\item $p = 3$, $C \cong \PP^1$, and the singular fibers of $f$ are of type $({\rm II},{\rm I}_{3^{2k}})$ or $({\rm II},{\rm I}^*_{3^{2k-1}})$ for some $k \geq 1$.
\item $p = 2$, $C \cong \PP^1$, and the singular fibers of $f$ are of type $({\rm II},{\rm I}_{2^{2k+1}})$ or $({\rm III},{\rm I}_{2^{2k+1}})$ for some $k \geq 1$.
\item $p \in \{2,3\}$, $C \cong \PP^1$, and $f$ is isotrivial with a unique singular fiber.
\end{enumerate}
\end{Theoremx}

One can classify the possible singular fibers of $f$ in Theorem (D) (ii). They are precisely the types described in Lemma \ref{twofibers} and Lemma \ref{twofibers23}.
All cases in Theorem (D) actually occur and we give the corresponding examples in Section \ref{sectionhorizontal}. Theorem (D) has the following well-known consequence (compare \cite[Theorem 7]{RudakovShafarevich}).
\begin{Corollary}\label{K3corollary}
There are no global regular vector fields on K3 surfaces.
\end{Corollary}

Our proof follows the strategy of Rudakov and Shafarevich and by using Theorem (D) we can circumvent the issues with Lemma 3 and Lemma 4 in \cite{RudakovShafarevich}. In characteristic $2$ and $3$, the proof builds on the fact that supersingular K3 surfaces admit an elliptic fibration with at least two singular fibers (see \cite[p.1502]{RudakovShafarevich2}). Let us also remark that there is an independent proof of Corollary \ref{K3corollary} due to Nygaard in \cite{Nygaard}.

Since our analysis of elliptic surfaces does not assume the existence of a section, we can also apply it to study the automorphism group scheme of Enriques surfaces in characteristic $2$. For example, we prove the following result in Example \ref{multiplicativeEnriques}.

\begin{Corollary}\label{Enriquescorollary}
Let $X$ be a generic supersingular Enriques surface in characteristic $2$. Then, $\Aut_X^0 \cong \mu_2$.
\end{Corollary}

Using the more refined results we give in Section \ref{vertical} and Section \ref{horizontal}, it is possible to determine the group scheme $\Aut_X^0$ in many concrete cases. For example, as an extension of Corollary \ref{Enriquescorollary}, we will use our results to calculate the connected components of the identity of the automorphism schemes of elliptic Enriques surfaces in characteristic $2$ in an upcoming article.

\begin{Remark}
All elliptic surfaces treated in this article are assumed to be relatively minimal. However, note that applying Blanchard's Lemma to the morphism $\pi: \widetilde{X} \to X$ from a relatively non-minimal elliptic surface $\tilde{f}: \widetilde{X} \to C$ to its relatively minimal model $f: X \to C$, we obtain an inclusion $\Aut_{\widetilde{X}}^0 \subseteq \Aut_X^0$.
In particular, suitably modified versions of Theorems (A), (B), (C), and (D) apply to $\widetilde{X}$ as well. We leave the formulation of these generalizations to the interested reader.
\end{Remark}

The outline of this article is as follows: In Section \ref{preliminaries}, we will give several general results on automorphism schemes of proper schemes, such as the behaviour under birational modifications, the relation to deformation theory, and a fixed point formula for actions of connected linearly reductive group schemes. This part of the article applies to arbitrary proper schemes over arbitrary fields and we hope that our results will help to answer Questions (A), (B), and (C) for more general classes of proper varieties. Then, we give some background on elliptic surfaces and recall the structure of the automorphism scheme of curves of genus one. In Section \ref{vertical}, we study the vertical component $\Ker(f_*)$ of $\Aut_X^0$ for an elliptic surface $X$ and in Section \ref{horizontal}, we describe the horizontal component $\im(f_*)$ of $\Aut_X^0$. Before deducing our main results from this in Section \ref{proofs}, we give several examples in Section \ref{examples}, illustrating the different phenomena that occur in the context of automorphism schemes of elliptic surfaces.

\vspace{3mm}
\noindent
\textbf{Acknowledgements:}
The author gratefully acknowledges funding from the DFG under research grant MA 8510/1-1.
Also, the author would like to thank the Department of Mathematics at the University of Utah for its hospitality while this article was written. I would like to thank Michel Brion for pointing me towards Blanchard's Lemma and Claudia Stadlmayr for fruitful discussions. Finally, I would like to thank Fabio Bernasconi, Michel Brion, Leo Herr, and Claudia Stadlmayr for helpful remarks on an earlier version of this article.

\section{Generalities} \label{preliminaries}

\subsection{Generalities on group scheme actions}

Throughout Section 2.1, we will be working over an arbitrary field $K$ and all schemes we consider will be of finite type over $K$. The following theorem of Matsumura and Oort \cite{MatsumuraOort} was mentioned in the introduction and proves the existence of our main object of interest.

\begin{Theorem} \label{MatsumuraOort}
Let $X$ be a proper scheme over $K$. Then, the functor
\begin{eqnarray*}
\Aut_{X/K}: &(Sch/K)^{op} &\to (Sets) \\
&S &\mapsto \Aut(X \times S \to S),
\end{eqnarray*}
where $\Aut(X \times S \to S)$ is the group of automorphisms of $X \times S$ over $S$, is representable by a group scheme $\Aut_{X/K}$ which is locally of finite type over $K$. In particular, its connected component of the identity $\Aut_X^0$ is of finite type over $K$.
\end{Theorem}

\begin{Remark}
More generally, one can prove the existence of a scheme of automorphisms $\Aut_{X/C}$ for a proper and flat morphism $f: X \to C$ of schemes, where $C$ is a normal and locally Noetherian scheme of dimension at most $1$, as follows: The relative Hilbert functor for $X \times_C X \to C$ is representable by a separated algebraic space $\mathcal{H}$ which is locally of finite presentation over $C$ (by \cite[Section 6]{ArtinI}), the functor  $\Aut_{X/C}$ is an open subfunctor of $\mathcal{H}$ (by \cite[Proposition 4.6.7 (ii)]{Ega3}), and every separated algebraic group space that is locally of finite type over $C$ is in fact a group scheme over $C$ (by \cite[Th\'eor\`eme (3.3.1).]{Raynaud}). This also shows that $\Aut_{X/C}$ exists as an algebraic group space under much weaker assumptions on $C$.
\end{Remark}

If the base scheme $C$ is clear from the context, we will simply write $\Aut_{X}$ for the functor of automorphisms of $X$ over $C$. For a scheme $X$ and a closed subscheme $Z \subseteq X$, we let $\Aut_{(Z,X)} \subseteq \Aut_X$ be the subgroup functor of automorphisms of $X$ preserving $Z$. Its $S$-valued points are given by
\begin{eqnarray*}
\Aut_{(Z,X)}(S) = \{ \alpha \in \Aut_X(S) \mid  Z \times S = (X \times S) \times_{\alpha,(X \times S)} (Z \times S)\}.
\end{eqnarray*}
Here, by $Z \times S = (X \times S) \times_{\alpha,(X \times S)} (Z \times S)$ we mean equality as closed subschemes of $X \times S$. Equivalently, $\Aut_{(Z,X)}$ is the stabilizer of the $K$-valued point corresponding to $Z$ in the Hilbert functor of $X$ over $K$. This second interpretation shows the following.

\begin{Lemma}
If $X$ is proper, then $\Aut_{(Z,X)}$ is a closed subgroup scheme of $\Aut_X$.
\end{Lemma}

\begin{Remark}\label{autzxremark}
If a group scheme $G$ acts on $X$, then the condition that $G \to \Aut_X$ factors through $\Aut_{(Z,X)}$ can be rephrased as $\rho^{-1}I_Z \cdot \mathcal{O}_{G \times X} = {\rm pr}_2^{-1}I_Z \cdot \mathcal{O}_{G \times X}$, where $\rho: G \times X \to X$ is the action, ${\rm pr}_2$ is the second projection, $I_Z$ is the ideal sheaf of $Z$ in $X$, and ${\rm pr}_2^{-1}I_Z \cdot \mathcal{O}_{G \times X}$ and $\rho^{-1}I_Z \cdot \mathcal{O}_{G \times X}$ denote the corresponding inverse image ideal sheaves.  For more details, see \cite[Section 2]{Fogarty}.
\end{Remark}
\vspace{1mm}

\subsubsection{Equivariant morphisms and birational modifications}

In this section, we will study how group scheme actions behave with respect to (birational) morphisms. Recall that the \emph{schematic image} of a morphism of schemes $f: X \to Y$ is the smallest closed subscheme of $Y$ through which $f$ factors. Since $X$ and $Y$ are of finite type over $K$, the schematic image can be described as the closed subscheme of $Y$ cut out by the quasi-coherent sheaf of ideals ${\rm Ker}(\mathcal{O}_Y \to f_* \mathcal{O}_X)$. The formation of $\Aut_{(Z,X)}$ is compatible with schematic images in the following sense.

\begin{Lemma}\label{image}
Let $Z \subseteq X$ be a closed subscheme of a scheme $X$ and $G$ a subgroup functor of $\Aut_{(Z,X)}$. Let $f: X \to Y$ be a $G$-equivariant morphism. Then, the morphism $G \to \Aut_Y$ factors through $\Aut_{(f(Z),Y)}$, where $f(Z)$ is the schematic image of $Z$ under $f$.
\end{Lemma}

\prf
It suffices to observe that for every $K$-scheme $S$ and $f_S: X \times S \to Y \times S$, the schematic images satisfy $f_S(Z \times S) = f(Z) \times S$. This is true by \cite[Lemma 14.6]{GoertzWedhorn}.
\qed
\vspace{3mm}

In general, not every $G$-action on $X$ descends to $Y$. An important case where we get an induced $G$-action is given by Blanchard's Lemma (see e.g. \cite[Theorem 7.2.1]{Brion}):

\begin{Theorem}[Blanchard's Lemma]\label{Blanchard}
Let $f: X \to Y$ be a proper morphism of schemes with $f_* \OO_X \cong \OO_Y$ and let $G$ be a connected group scheme acting on $X$. Then, the following hold:
\begin{enumerate}[(i)]
\item There is a unique $G$-action on $Y$ such that $f$ is $G$-equivariant.
\item If $X$ and $Y$ are proper, there is a natural homomorphism $f_*: \Aut_X^0 \to \Aut_Y^0$.
\item If, additionally, $X$ and $Y$ are integral and $f$ is birational, then $f_*: \Aut_X^0 \to \Aut_Y^0$ is a closed immersion.
\end{enumerate}
\end{Theorem}

Alternatively, we can start with a group scheme action on $Y$ and ask whether it lifts along a birational modification $f: X \to Y$ to a compatible action on $X$.

\begin{Proposition} \label{blowup}
Let $Z \subseteq X$ be a closed subscheme of a scheme $X$ and let $\pi: {\rm Bl}_Z(X) \to X$ be the blow-up of $X$ in $Z$. Let $G$ be a group scheme acting on $X$. If $G \to \Aut_X$ factors through $\Aut_{(Z,X)}$, then the $G$-action lifts to ${\rm Bl}_Z(X)$. The converse holds if $G \to \Aut_{{\rm Bl}_Z(X)}$ factors through $\Aut_{(\pi^{-1}(Z),{\rm Bl}_Z(X))}$, $Z$ contains no irreducible component of $X_{red}$, and either $Z$ is reduced or $Z \subseteq X$ is a regular embedding.
\end{Proposition}

\prf
Since the action map $\rho: G \times X \to X$ is flat and blow-up commutes with flat base-change, we have the following diagram of solid arrows with cartesian square
$$
\xymatrix{
 G \times {\rm Bl}_Z(X) \ar[dr]^{{\rm id} \times \pi} \ar@{.>}[r]^(.32){\iota} & Y := {\rm Proj} \bigoplus\limits_{i = 0}^{\infty} (\rho^{-1}I_Z \cdot\mathcal{O}_{G \times X})^i \ar[r]^(.72){\rho'} \ar[d]^{\pi'} &  {\rm Bl}_Z(X) \ar[d]^\pi \\
 & G \times X \ar[r]^\rho & X
}
$$

\noindent and we are asking for the existence of the map $\iota$ such that $\rho' \circ \iota$ is an action of $G$ and such that the above diagram commutes. If $G \to \Aut_X$ factors through $\Aut_{(Z,X)}$, then $Y \cong G \times {\rm Bl}_Z(X)$ over $G \times X$ by Remark \ref{autzxremark} and we get the desired map $\iota$. 

For the converse, we apply Lemma \ref{image}, which shows that it suffices to check that $Z$ is the schematic image of $\pi^{-1}(Z)$. 
Recall that $\pi^{-1}(Z) = \Proj \bigoplus_{i=0}^{\infty} (I_Z^{i}/I_{Z}^{i+1})$. If $Z$ does not contain an irreducible component of $X_{red}$, then $I_Z$ is not nilpotent at any $z \in Z$, hence the fibers of $\pi|_{\pi^{-1}(Z)}: \pi^{-1}(Z) \to Z$ are non-empty. If $Z$ is reduced, this implies that $Z$ is the schematic image of $\pi^{-1}(Z)$. If $Z \subseteq X$ is a regular embedding, then $\pi^{-1}(Z) \to Z$ is a projective bundle, hence also in this case $Z$ is the schematic image of $\pi^{-1}(Z)$.
\qed
\vspace{1mm}

\begin{Remark}
In Proposition \ref{criteria}, we will see that if $G$ is connected and the normal bundle $N_{E/{\rm Bl}_Z(X)}$ of $E := \pi^{-1}(Z)$ in ${\rm Bl}_Z(X)$ satisfies $h^0(E,N_{E/{\rm Bl}_Z(X)}) = 0$, then $G$ preserves $E$. In particular, this holds if $\pi$ is the contraction of a negative definite configuration of curves on a smooth surface.
\end{Remark}

Let $\nu: \widetilde{X} \to X$ be a finite and birational morphism. Then, the \emph{conductor ideal} $I_\nu$ of $\nu$ is defined as 
$$
I_\nu = \Hom_{\OO_X}(\nu_* \OO_{\widetilde X}, \OO_X)  = {\rm Ann}_{\OO_X}(\nu_* \OO_{\widetilde X}/\OO_X)\subseteq \OO_X.
$$
We let $C_\nu \subseteq X$ be the closed subscheme defined by $I_\nu$ and call it the \emph{conductor} of $\nu$. Then, the locus where $\nu$ is not an isomorphism is precisely $C_\nu$.
If both $\widetilde{X}$ and $X$ are Gorenstein, it follows from relative duality \cite{Kleiman} that $I_\nu \cdot \omega_X = \nu_* (\omega_{\widetilde{X}})$. In particular, $I_\nu$ is reflexive and thus $C_\nu$ is a generalized divisor in the sense of \cite{Hartshornegeneralized}. The fact that $\nu^{-1} I_\nu \cdot \mathcal{O}_{\widetilde{X}}$ is locally principal can be used to show that $\nu$ is the blow-up of $I_\nu$ (see \cite[Proposition 2.9]{Piene}). Using Proposition \ref{blowup}, we obtain the following proposition.

\begin{Proposition}\label{normalization}
Let $\nu: \widetilde{X} \to X$ be a finite and birational morphism between Gorenstein schemes. Then, $\widetilde{X} = {\rm Bl}_{C_\nu}(X)$ and thus a $G$-action on $X$ lifts to $\widetilde{X}$ if and only if $G \to \Aut_X$ factors through $\Aut_{(C_\nu,X)}$.
\end{Proposition}

\prf
By Proposition \ref{blowup}, the $G$-action on $X$ lifts to  $\widetilde{X}$ if $G \to \Aut_X$ factors through $\Aut_{(C_\nu,X)}$.

For the converse, note that, because $\rho$ and ${\rm pr_2}: G \times X \to X$ are flat and ${\rm id} \times \nu$ is the base change of $\nu$ along both $\rho$ and ${\rm pr_2}$, we can use the fact that cohomology, annihilators, and quotients commute with flat base change to obtain
\begin{eqnarray*}
{\rm pr_2}^{-1}(I_{C_\nu}) \cdot\mathcal{O}_{G \times X} &=& {\rm pr_2}^{-1}({\rm Ann}_{\OO_X}(\nu_* \OO_{\widetilde X}/\OO_X)) \cdot\mathcal{O}_{G \times X}
\\ &=&  {\rm Ann}_{\OO_{G \times X}}((({\rm id} \times \nu)_*\OO_{G \times {\widetilde X}})/\OO_{G \times X}) \cdot\mathcal{O}_{G \times X} \\
 &=& \rho^{-1}({\rm Ann}_{\OO_X}(\nu_* \OO_{\widetilde X}/\OO_X)) \cdot\mathcal{O}_{G \times X} = \rho^{-1}(I_{C_\nu}) \cdot\mathcal{O}_{G \times X}.
\end{eqnarray*}
Hence, by Remark \ref{autzxremark}, the $G$-action on $X$ preserves $C_{\nu}$.
\qed

\begin{Remark}
In particular, if $X$ is a reduced proper scheme over $K$ such that $X$ and its normalization $\widetilde{X}$ are Gorenstein, the scheme of automorphisms of $X$ that lift to $\widetilde{X}$ is precisely the stabilizer of the conductor. This seems to be the "general principle" mentioned in the calculation of the automorphism scheme of a cuspidal plane cubic curve in \cite[p. 213]{BombieriMumford3}.
\end{Remark}
\vspace{0.5mm}

\subsubsection{Fixed points}
Recall that if a group scheme $G$ acts on a scheme $X$, then the subfunctor of fixed points for this action is defined as
$$
X^G(S) := \{ x \in X(S) \mid g(x_T) = x_T \text{ for all }S\text{-schemes } T \text{ and } g \in G(T)\}
$$
By \cite[Theorem 2.3]{Fogarty}, $X^G$ is representable by a closed subscheme of $X$. We have the following lemma, whose proof is the same as the one of the corresponding statement for actions of abstract groups and thus left to the reader.

\begin{Lemma} \label{fixedofsubscheme}
Let $X$ be a scheme and let $H \subseteq G$ be a normal subgroup scheme of a group scheme $G$. Assume that $G$ acts on $X$. Then, $G$ acts on the fixed locus of the induced $H$-action on $X$, that is, $G \to \Aut_X$ factors through $\Aut_{(X^H,X)}$.
\end{Lemma}

This simple observation can sometimes be used to obtain information about fixed points of $G$ via the following proposition.

\begin{Proposition} \label{isolatedefixedpoints}
Let $X$ be a scheme and let $H \subseteq G$ be a normal subgroup scheme of a connected group scheme $G$ that acts on $X$. Assume that $X^H$ admits a connected component $P$ isomorphic to $\Spec K$. Then, $P \in X^G$
\end{Proposition}

\prf
By Lemma \ref{fixedofsubscheme}, $G$ acts on $X^H$ and since $G$ is connected, this action preserves the connected components of $X^H$. By our assumption, the connected component of $X^H$ containing $P$ is isomorphic to $\Spec K$. Therefore, the induced $G$-action on $P$ is trivial, hence $P \in X^G$.
\qed

\vspace{5mm}

\subsubsection{Some deformation theory}
In this section, we will use the deformation theory of a closed subscheme $Z$ of a scheme $X$ to obtain information about the functor ${\rm Aut}_{(Z,X)}$.
For the necessary background on deformation theory, we refer the reader to \cite{Sernesi}. 
We fix the following notation.

\begin{Notation}\label{Deformationnotation}
Let $Z \subseteq X$ be a closed subscheme of a scheme $X$.
\begin{itemize}
  \item $\Def_Z$ is the functor of deformations of $Z$.
  \item ${\rm Def}_{Z/X}$ is the functor of deformations of $Z$ in $X$ for which the deformation of the ambient scheme $X$ is trivial.
  \item ${\rm Def}'_{Z/X}$ is the subfunctor of deformations of $Z$ in $X$ as above mapping to the trivial deformation of $Z$ via the forgetful map $F: {\rm Def}_{Z/X} \to \Def_Z$.
   \item $\widehat{\Aut}_{X}$ is the restriction of $\Aut_X$ to the category ${\rm Art}_K^{\rm op}$ of Artinian local $K$-schemes with residue field $K$ whose closed points map to ${\rm id}_X$. For every such $S \in {\rm Art}_K^{\rm op}$, there is a natural map $\widehat{\Aut}_{X}(S) \to {\rm Def}'_{Z/X}(S)$ given by $\alpha \mapsto (X \times S) \times_{\alpha,(X \times S)} (Z \times S)$.
\end{itemize}
\end{Notation}

Recall that a functor of Artin rings $F: {\rm Art}_K \to ({\rm Sets})$ is said to be \emph{prorepresentable} by a complete Noetherian local $K$-algebra $R$ with residue field $K$ if there exists a natural isomorphism ${\rm Hom}_K(R,-) \to F$. A \emph{hull} for $F$ is an $R$ as above together with a formally smooth natural transformation ${\rm Hom}_K(R,-) \to F$ that induces an isomorphism on tangent spaces.

\begin{Lemma}
Let $Z \subseteq X$ be a closed subscheme of a proper scheme $X$. Then, $\Def_X$ has a hull and $\Def_{Z/X}$ is prorepresentable. If $H^0(X,T_X) = 0$, then $\Def_X$ is prorepresentable.
\end{Lemma}

The situation for $\Def'_{Z/X}$ is more subtle. Using Schlessinger's criteria \cite{Schlessinger}, one can prove that $\Def'_{Z/X}$ is prorepresentable if and only if it has a hull. But even if $X$ is smooth and proper, the functor $\Def'_{Z/X}$ is not prorepresentable in general. In geometric terms, this means that for a deformation of $Z$ in $X$ over an Artinian local $K$-scheme $S$, the locus over which the deformation of $Z$ is trivial once we forget about its embedding into $X_S$ need not be a scheme.
Such a phenomenon can only occur if the deformation theory of $Z$ itself is pathological: If $H$ is a hull for $\Def_Z$ and $\{\ast\} \to \Def_Z$ is the morphism that maps $S$ to the trivial deformation, then the base change $H' := \{\ast\} \times_{\Def_Z} H \to H$ is a monomorphism of functors of Artin rings but $H'$ may not be prorepresentable. If $H'$ is prorepresentable, that is, if the locus in $H$, over which the deformation of $Z$ is \underline{t}rivial, is \underline{r}epresentable, we say that $H$ is a \emph{TR-hull} for $\Def_Z$.
Note that if $\Def_Z$ is prorepresentable, then it admits a TR-hull.

\begin{Lemma}\label{prorep}
Let $Z \subseteq X$ be a closed subscheme of a proper scheme $X$. If $\Def_Z$ admits a TR-hull, then $\Def_{Z/X}'$ is prorepresentable. 
\end{Lemma}

\prf
By definition, $\Def_{Z/X}'$ is the fiber product $\{\ast\} \times_{\Def_Z} \Def_{Z/X}$. If $\Def_Z$ admits a TR-hull $H$ with $H' := \{\ast\} \times_{\Def_Z} H$, then $\Def_{Z/X}' = H' \times_H \Def_{Z/X}$ and thus $\Def_{Z/X}'$ is prorepresentable.
\qed

\begin{Example} \label{cuspexample}
If $Z$ is a cuspidal plane cubic over an algebraically closed field $k$ of characteristic $p \geq 0$, then $\Def_Z$ admits a TR-hull if and only if $p \not\in \{2,3\}$ (see \cite[p.202]{BombieriMumford3}).
\end{Example}

The reason why we care about the functor $\Def_{Z/X}'$ is that we can use it to check whether the inclusion $\Aut_{(Z,X)}^0 \subseteq \Aut_X^0$ is an equality.

\begin{Lemma} \label{restriction}
Let $X$ be a proper scheme and $Z \subseteq X$ a closed subscheme. The natural map $\widehat{\Aut}_{X} \to {\rm Def}'_{Z/X}$ is constant if and only if $\Aut_{(Z,X)}^0 = \Aut_X^0$.
\end{Lemma}

\prf
Since $\Aut_{(Z,X)}^0 \to \Aut_X^0$ is a closed immersion and both sides are connected, the equality $\Aut_{(Z,X)}^0 = \Aut_X^0$ holds if and only if $\widehat{\Aut}_{X} = \widehat{\Aut}_{(Z,X)}$.
From the definitions, we see that $\widehat{\Aut}_{(Z,X)}$ is the fiber of $\widehat{\Aut}_{X} \to {\rm Def}'_{Z/X}$ over the trivial deformation of $Z$ in $X$. Thus, $\widehat{\Aut}_{X} = \widehat{\Aut}_{(Z,X)}$ if and only if $\widehat{\Aut}_{X} \to {\rm Def}'_{Z/X}$ is constant.
\qed
\vspace{5mm}

Now, we want to understand the tangent spaces of the functors recalled in Notation \ref{Deformationnotation}. To this end, we define a subsheaf $T_X\langle Z \rangle \subset T_X$ of the tangent sheaf $T_X$ of $X$ via
$$
T_X\langle Z \rangle(U) = \{ D \in T_X(U) \mid D(I_Z(U)) \subseteq I_Z(U) \}.
$$ 
We recall that $N_{Z/X} := (I_Z/I_Z^2)^{\vee}$ denotes the normal sheaf of $Z$ in $X$ and that $K[\epsilon] := K[x]/x^2$ is the ring of dual numbers. Then, with a slight abuse of notation, we get the following well-known identifications of the relevant tangent spaces:
\begin{itemize}
 \item $\widehat{\Aut}_{(Z,X)}(K[\epsilon])  = H^0(X,T_X\langle Z \rangle)$.
 \item $\widehat{\Aut}_{X}(K[\epsilon]) = H^0(X,T_X)$.
 \item $\Def_{Z/X}(K[\epsilon]) = H^0(Z,N_{Z/X})$.
 \item $\Def'_{Z/X}(K[\epsilon]) = \ker(H^0(Z,N_{Z/X}) \to \Def_Z(k[\epsilon]))$.
 \item If $Z$ is reduced, then $\Def_Z(K[\epsilon]) = {\rm Ext}^1(\Omega_Z,\OO_Z)$ and if $Z$ is smooth, then ${\rm Ext}^1(\Omega_Z,\OO_Z) = H^1(Z,T_Z)$. In these cases, the differential of the forgetful map $F: \Def_{Z/X} \to \Def_Z$ is induced by the conormal sequence.
\end{itemize}

\begin{Proposition} \label{criteria}
 If $\Def'_{Z/X}$ is trivial and $X$ is proper, then $\Aut_{(Z,X)}^0 = \Aut_X^0$. This holds in each of the following cases:
\begin{enumerate}[(a)]
 \item $H^0(Z,N_{Z/X}) = 0$.
 \item $Z$ is reduced, $\Def_{Z/X}'$ is prorepresentable and $H^0(Z,N_{Z/X}) \to {\rm Ext}^1(\Omega_Z,\OO_Z)$ is injective.
 \item $X$ is smooth in a neighborhood of $Z$ and $Z$ is a geometrically reduced, geometrically connected, and non-smooth effective divisor on $X$ with $N_{Z/X} = \OO_Z$ such that $\Def_{Z/X}'$ is prorepresentable.
\end{enumerate}

\end{Proposition}

\prf
First, observe that Lemma \ref{restriction} shows that  $\Aut_{(Z,X)}^0 = \Aut_X^0$ holds if $\Def'_{Z/X}$ is trivial, so we have to check that $\Def'_{Z/X}$ is trivial under any of the stated conditions.

Since $\Def_{Z/X}$ is prorepresentable, it is trivial as soon as $H^0(Z,N_{Z/X}) = 0$. As $\Def'_{Z/X}$ is a subfunctor of $\Def_{Z/X}$, it is also trivial in this case. This is Claim (a).

As for Claim (b), since $\Def'_{Z/X}$ is prorepresentable, it suffices to check that $\Def'_{Z/X}(K[\epsilon])$ is trivial. But by the facts recalled above and since $Z$ is reduced, we have $\Def'_{Z/X}(K[\epsilon]) = \ker(H^0(Z,N_{Z/X}) \to {\rm Ext}^1(\Omega_Z,\OO_Z))$.

To prove Claim (c), we thus have to prove that $H^0(Z,N_{Z/X}) \to {\rm Ext}^1(\Omega_Z,\OO_Z)$ is injective. Since $Z$ is a reduced effective Cartier divisor, we have the short exact conormal sequence
$$
0 \to \OO_Z(Z) \to \Omega_X|_Z \to \Omega_Z \to 0.
$$
Applying $\Hom(-,\OO_Z)$, we obtain 
$$
H^0(Z,T_X|_Z) \to H^0(Z,N_{Z/X}) = \Hom(\OO_Z(Z),\OO_Z) \overset{f}{\to} \Ext^1(\Omega_Z,\OO_Z).
$$
The map $f$ associates to a morphism $\varphi: \OO_Z(Z) \to \OO_Z$ the pushout of the conormal sequence along $\varphi$. Since $Z$ is geometrically reduced and geometrically connected, the space $H^0(Z,N_{Z/X}) = H^0(Z,\OO_Z)$ is $1$-dimensional, so that $f$ is either trivial or injective. 
Suppose that $f$ is trivial, that is, that $f({\rm id}) = 0$. This means that the conormal sequence splits. Thus, $\Omega_Z$ is locally free, being a direct summand of the locally free sheaf $\Omega_X|_Z$ and the rank of $\Omega_Z$ is $\dim Z = \dim X - 1$. Therefore, $Z$ is smooth, contradicting our assumption that $Z$ is non-smooth. Hence, $f$ is injective and Claim (c) follows from Claim (b).
\qed
\vspace{3mm}

\begin{Remark} \label{criteriaremark}
Without the assumption on the prorepresentabilty of $\Def_{Z/X}'$, the proof of Proposition \ref{criteria} (b) and (c) shows that $H^0(X,T_X\langle Z \rangle) = H^0(X,T_X)$. Indeed, the map $\varphi: H^0(X,T_X) \to H^0(Z,N_{Z/X})$ factors through $H^0(Z,T_X|_Z)$, so the above proof shows that $\varphi$ is trivial, hence $H^0(X,T_X\langle Z \rangle) = \Ker(\varphi) = H^0(X,T_X)$. In particular, even if $\Def_{Z/X}'$ is not prorepresentable, the functors $\Aut_{(Z,X)}$ and $\Aut_{X}$ have the same tangent space at the identity in case (b) and (c). This implies, for example, that for every connected subgroup scheme $G$ of $\Aut_X$ the intersection $G \cap \Aut_{(Z,X)}$ is non-trivial.
\end{Remark}

Another case where $\Aut_{(Z,X)}$ and $\Aut_X$ have the same tangent space is if $Z$ is given by a Frobenius power of an ideal. Recall that if $I \subseteq \OO_X$ is an ideal sheaf, then its \emph{Frobenius power} $I^{[p]}$ is the ideal sheaf which is locally generated by the $p$-th powers of generators of $I$. If $I$ is locally principal, then $I^{[p]} = I^p$.

\begin{Lemma} \label{Frobeniuspower}
Let $Z \subseteq X$ be a closed subscheme of a scheme $X$ and let $Z^{[p]} \subseteq X$ be the closed subscheme defined by $I_Z^{[p]}$. Then,  $T_X\langle Z^{[p]} \rangle = T_X$. In particular, $H^0(X,T_X\langle Z^{[p]} \rangle) = H^0(X,T_X)$ holds.
\end{Lemma}

\prf
Let $U \subseteq X$ be an open subset with $I_Z(U) = \langle f_1,\hdots,f_n \rangle$ and $D \in T_X(U)$. Then, using the Leibniz rule, we deduce for arbitrary $a_i \in \OO_X(U)$ that
$$
D(\sum_{i=1}^n a_i f_i^p) = \sum_{i=1}^n f_i^p D(a_i) \in I_{Z^{[p]}}(U).
$$
\vspace{-5mm}
\qed
\vspace{5mm}

\subsubsection{Examples of group schemes and some structure theory}
The material in this section is standard and we refer the reader to \cite{SGA3}, \cite{SGA3II}, \cite{Brion}, \cite{Milne}, and \cite{Waterhouse} for proofs. 

If $\Char(K) > 0$ and $X$ is a scheme over $K$, we write $X^{(p)}$ for the pullback of $X$ along the $K$-linear Frobenius. For a group scheme $G$ over a field of positive characteristic $K$, the notation $G[F^n]$ denotes the kernel of the $n$-fold $K$-linear Frobenius $F^n: G \to G^{(p^n)}$. If $G$ is finite and connected, then $G[F^n] = G$ for $n \gg 0$. If $k \subseteq K$ is a field extension, we write $G_K$ for $G \times_{\Spec k} \Spec K$.  Let us recall $G[F^n]$ for the smooth $1$-dimensional group schemes over an algebraically closed field $k$:

\begin{Lemma}
Let $k$ be an algebraically closed field. Let $G$ be a smooth $1$-dimensional group scheme over $k$. Then, $G[F^n]$ is as follows:
\begin{itemize}
\item If $G = \GG_m$, then $G[F^n] = \mu_{p^n}$.
\item If $G = \GG_a$, then $G[F^n] = \alpha_{p^n}$.
\item If $E$ is an \emph{ordinary} elliptic curve, then $E[F^n] = \mu_{p^n}$.
\item If $E$ a \emph{supersingular} elliptic curve over $k$, then $E[F^n] = M_n$, where $M_n$ is an $n$-fold non-split extension of $\alpha_p$ by itself.
\end{itemize}
\end{Lemma}

In each case, ${\rm length}(G[F^n]/G[F^{n-1}]) = p$, so the above list is a complete enumeration of all finite connected subgroup schemes of these four group schemes. 

Next, let us recall some structural results for a group scheme $G$ of finite type over a field $K$.

\begin{Lemma}\label{generalgroupscheme}
Let $G$ be a group scheme of finite type over a field $K$. Then,
\begin{enumerate}[(i)]
\item \emph{(Cartier's Theorem)} If $\Char(K) = 0$, then $G$ is reduced. 
\item The connected component of the identity $G^0 \subseteq G$ is a closed subgroup scheme.
\item There is a smallest normal subgroup scheme $H \subseteq G$ such that $G^{ab} := G/H$ is commutative.
\item If $K$ is perfect, then $G_{red} \subseteq G$ is a closed and smooth subgroup scheme.
\end{enumerate}
\end{Lemma}

Denoting by $\mathcal{H}om(G,H)$ the sheaf of homomorphisms between two group schemes $G$ and $H$, recall that the \emph{Cartier dual} $G^\vee := \mathcal{H}om(G,\GG_m)$ of a finite commutative group scheme $G$ is also a finite and commutative group scheme and we have $(G^\vee)^\vee \cong G$. Moreover, for abelian group schemes, recall that there exists a morphism $V: G^{(p)} \to G$ called \emph{Verschiebung}, which, if $G$ is finite, coincides with the morphism induced by the Frobenius on $G^\vee$.

\begin{Lemma}\label{finitecommgroupscheme}
Let $G$ be a finite and commutative group scheme over a field $K$. Then, the following hold:
\begin{enumerate}[(i)]
\item If $K$ is perfect, then there is a functorial decomposition
$$
G \cong G_{rr} \times G_{rl} \times G_{lr} \times G_{ll}
$$
where $G_{xy}$ is reduced if $x = r$ and connected if $x = l$, and $G_{xy}^\vee$ is reduced if $y = r$ and connected if $y = l$. We say that $G$ is of type $xy$ if $G \cong G_{xy}$ for $x,y \in \{r,l\}$.
\item If $K = k $ is algebraically closed, then 
\begin{enumerate}[(1)]
\item $G_{rr}$ is the constant group scheme associated to an abelian group of order prime to $p$,
\item $G_{rl}$ is the constant group scheme associated to an abelian group of $p$-power order,
\item $G_{lr} \cong \prod_{i=1}^m \mu_{p^{n_i}}$ for some $n_i,m \geq 0$,
\item $G_{ll}$ is an iterated extension of $\alpha_p$ by itself. Moreover, $G_{ll} \cong \alpha_p^r$ for some $r \geq 0$ if and only if both $F$ and $V$ are trivial on $G_{ll}$.
\end{enumerate}
\end{enumerate}
\end{Lemma}

\begin{Lemma}\label{homsandexts}
Let $G$ and $H$ be finite and commutative group schemes over an algebraically closed field $k$. Then, the following hold:
\begin{enumerate}[(i)]
\item If $G$ is of type $xy$ and $H$ is of type $x'y'$ and if there is a non-trivial homomorphism $G \to H$, then $(x,y) = (x',y')$.
\item If $M$ is an extension of $G$ by $H$, then we have the following:
\begin{enumerate}[(1)]
\item If $G$ and $H$ are of type $lr$, then $M$ is commutative of type $lr$.
\item If $G = \alpha_p$ and $H$ is of type $lr$, then $M \cong H \times G$.
\item If $G$ is of type $lr$ and $H$ is of type $ll$, then $M \cong H \rtimes G$.
\end{enumerate}
\end{enumerate}
\end{Lemma}

\prf
Claim (i) follows from functoriality of the canonical decomposition of a finite commutative group scheme. Claim (ii) (1) is \cite[Theorem 15.39.]{Milne}, the splitting in Claim (ii) (2) follows from \cite[Th\'eor\`eme 6.1.1. B), Exp. XVII]{SGA3II} and that $M$ is in fact a direct product follows from the fact that $\Aut_H$ is \'etale. Finally, Claim (ii) (3) follows from \cite[Th\'eor\`eme 5.1.1., Exp. XVII]{SGA3II}.
\qed
\vspace{5mm}

\subsection{Linearization of $\mu_{p^n}$-actions and a fixed point formula}\label{Linearization} From now on, we will work over an algebraically closed field $k$ of characteristic $p > 0$.
Let $X$ be a smooth variety over $k$ with a faithful $\mu_{p^n}$-action. Let $P \in X$ be a fixed point of this action. Since $\mu_{p^n}$ is linearly reductive, it is well-known (see e.g. \cite[Proof of Corollary 1.8]{Satriano}) that the action of $\mu_{p^n}$ on $X$ can be linearized in a formal neighborhood of $P$ in $X$. If $X$ is a surface, "linearizability" means that there is a $\mu_{p^n}$-equivariant isomorphism
$$
\widehat{\OO}_{X,P} \to k[[x,y]],
$$
where the action of $\mu_{p^n}$ on $k[[x,y]]$ is linear. After possibly conjugating this action, we can assume that $\mu_{p^n}$ acts diagonally on $k[[x,y]]$, that is, via the coaction
\begin{eqnarray*}
k[[x,y]] &\to& k[[x,y]] \otimes_k k[\lambda]/(\lambda^p -1) \\
x & \mapsto & x \otimes \lambda \\
y & \mapsto & y \otimes \lambda^a
\end{eqnarray*}
for some $0 \leq a \leq p^n-1$. We say that the $\mu_{p^n}$-action on $X$ is of type $\frac{1}{p^n}(1,a)$ at $P$. The fact that one can linearize $\mu_{p^n}$-actions has the following consequences for the fixed locus $X^{\mu_{p^n}}$  and the quotient $q: X \to Y:= X/\mu_{p^n}$.

\begin{itemize}
\item $X^{\mu_{p^n}}$ is representable by a smooth closed subscheme of $X$ (see \cite[Theorem 5.4]{Fogarty}).
\item If $X$ is a surface, then $q(P)$ is a singular point of $Y$ if and only if $P$ is an isolated fixed point of the $\mu_{p^n}$-action (see \cite[Theorem 1.3]{Satriano}).
\item If $X$ is a surface and $p \nmid a$, then $q(P)$ is a cyclic quotient singularity of type $\frac{1}{p^n}(1,a)$ (see e.g. \cite[Theorem 2.3]{Hirokado} for $n = 1$; the general case is similar).
\end{itemize}

\begin{Remark}\label{vectorfield}
Actions of $\mu_p$ (and $\alpha_p$) on $X$ correspond bijectively to vector fields $D$ on $X$ with $D^p = D$ (resp. $D^p = 0$). An explicit description of this correspondence can be found for example in \cite[Section 3]{Tziolas3}. We remark that the fixed locus of the action is identified with the zero locus of the vector field via this correspondence. Vector fields with $D^p = \lambda D$ for some $\lambda \in k(X)$ are called \emph{$p$-closed}, and $D$ is called \emph{multiplicative} (resp. \emph{additive}) if $D^p = D$ (resp. $D^p = 0$).
\end{Remark}
We will now prove a fixed point formula for $\mu_{p^n}$-actions on smooth projective varieties. It may be possible to give a proof similar to the proof of the fixed point formula for torus actions on smooth varieties by Iversen \cite{Iversen}, but we were not able to find a suitable reference.

\begin{Theorem} \label{fixedpointformula}
Let $X$ be a smooth proper variety with an action of a finite commutative group scheme $G$ of type $lr$. Then, the $\ell$-adic Euler characteristics of $X$ and $X^{G}$ coincide, that is, 
$$
e(X) = e(X^{G})
$$
\end{Theorem}

\prf
Since $X^{\mu_{p^m}}$ is a smooth proper variety for all $\mu_{p^m} \subseteq G$ and $G$ is commutative, we can use Lemma \ref{fixedofsubscheme}, Lemma \ref{finitecommgroupscheme} (ii) (3), and induction on the length of $G$ to assume without loss of generality that $G \cong \mu_p$.

Then, by Remark \ref{vectorfield}, the $G$-action corresponds to the action of a multiplicative vector field $D$, and the fixed locus $X^{\mu_p}$ coincides with the zero locus of $D$. Phrased differently, the fixed locus $X^{\mu_p}$ is the zero locus of a section of $T_X$ and therefore $e(X^{\mu_p}) = c_{\dim(X)}(T_X) = e(X)$.
\qed

\begin{Remark}
The fixed point formula for torus actions given in \cite{Iversen} can be deduced from Theorem \ref{fixedpointformula} by considering the $\mu_{p^n}$-actions induced by a given $\GG_m$-action and letting $n \to \infty$.
\end{Remark}

\begin{Remark}\label{alphapfixedpoint}
Note that if $X$ is a smooth proper variety with an action of $\alpha_p$, then also $e(X) = e(X^{\alpha_p})$, by the same proof as in the $\mu_p$-case. However, as $X^{\alpha_p}$ is not necessarily smooth, it is unclear how to extend this to actions of, say, $M_n$.
\end{Remark}

\vspace{5mm}

\subsection{Elliptic surfaces}
In this section, we will recall the necessary background on elliptic surfaces over the algebraically closed field $k$.
Following \cite{Mumford}, we say that a non-zero effective divisor $F = \sum_{i=1}^n a_iF_i$ on a smooth surface $X$ is \emph{of canonical type} if $F.F_i = K_X.F_i = 0$ for all $i$. We say that $F$ is \emph{indecomposable}, if it is not a non-trivial sum of divisors of canonical type. Every fiber of an elliptic fibration is a curve of canonical type and, conversely, for many surfaces, curves of canonical type can be used to prove the existence of elliptic fibrations.

Let us recall the Kodaira--N\'eron classification of indecomposable divisors $F$ of canonical type (see e.g. \cite{KodairaII}). If $F$ is irreducible, it is either an elliptic curve (Type ${\rm I}_0$), a nodal rational curve (Type ${\rm I}_1$) or a cuspidal rational curve (Type ${\rm II}$). If $F$ is not irreducible, its components are $(-2)$-curves. If the components of $F$ do not intersect transversally, then $F$ consists either of two $(-2)$-curves which meet with multiplicity $2$ at one point (Type ${\rm III}$) or of three $(-2)$-curves meeting transversally in a single point (Type ${\rm IV}$). In all other cases, all curves intersect transversally in distinct points and the resulting dual graphs are given in the following table. We call $F$ \emph{additive} if it is not of type ${\rm I}_n$ and \emph{multiplicative} if it is of type ${\rm I}_n$ with $n \geq 1$.

\vspace{-5mm}
$$
\begin{array}[t]{|l*{10}{|c}|}
\hline
\text{Type of $F$}             
& {\rm I}_0 & {\rm I}_1 &  {\rm I}_m    & \p{\rm I}_m^*    & {\rm II}    & {\rm III} & {\rm IV} &  \p{\rm IV}^*    & \p{\rm III}^* & \p{\rm II}^* \\
 \hline&&&&&&&&&&\\[-2ex]
\text{Dual Graph} & - & - &  \widetilde{A}_{m-1}    & \widetilde{D}_{m+4}    & - & \widetilde{A}_1 & \widetilde{A}_2 &  \widetilde{E}_6    & \widetilde{E}_7 & \widetilde{E}_8 \\
 \hline
\end{array}
$$
\vspace{2mm}

Recall that by a result of Lichtenbaum and Shafarevich, the minimal proper regular model of a curve of positive genus over a Dedekind scheme exists and is unique (see \cite[Theorem 9.3.21]{Liu}). In the setting of elliptic surfaces, this can be rephrased as follows.

\begin{Lemma}\label{uniquemodel}
Let $C$ be a smooth projective curve over $k$ and let $F_\eta \to \Spec k(C)$ be a smooth projective curve of genus $1$ over $k(C)$. Then, there exists a unique elliptic surface $f: X \to C$ with generic fiber $F_\eta$.
\end{Lemma}
Using this, the \emph{Jacobian} $J(f):J(X) \to C$ of an elliptic surface $f:X \to C$ is simply the minimal proper regular model of the Jacobian $\Pic^0_{F_\eta}$ of the generic fiber $F_\eta$ of $f$. Note that the line bundle $\OO_{F_\eta}$ induces a canonical section of $J(f)$ and, away from the multiple fibers of $f$, the smooth locus of $f$ is a torsor under the smooth locus of its Jacobian. We call an elliptic surface $f$ \emph{Jacobian} if $f$ admits a section.

Now, let $f: X \to C$ be an elliptic surface and choose integers $m_i$ and indecomposable divisors $F_i$ of canonical type for $i=1,\hdots,n$ such that the $m_iF_i$ are precisely the multiple fibers of $f$. Then, $m_i$ is called \emph{multiplicity} of $F_i$ and whenever we say that $mF$ is a fiber of $f$, we implicitly assume that $m$ is the multiplicity of the fiber.
 Set $P_i := f(F_i)$. 
Let $\cal{L} \oplus \cal{T}$ be the decomposition of $R^1f_* \OO_X$ into its locally free part $\cal{L}$ and its torsion part $\cal{T}$.  A multiple fiber $m_iF_i$ of $f$ is called \emph{wild} if $\cal{T}_{P_i} \neq 0$ and \emph{tame} otherwise. Equivalently, $m_iF_i$ is tame if and only if $\nu_i = m_i$, where $\nu_i$ is the order of the normal bundle of $F_i$ in $X$.

Recall the following formulas:

\begin {itemize}
 \item (\emph{Application of Riemann--Roch})
 $$\chi(X,\OO_X) = \chi(X, \omega_X^{\otimes n}) \geq 0$$ for all $n \in \mathbb{Z}$.
 \item (\emph{Noether formula})
$$12 \chi(X,\OO_X) = c_2(X).$$
 \item (\emph{Igusa inequality})
$${\rm rk}(\Num(X)) \leq b_2(X)$$.
 \item (\emph{Canonical bundle formula} (see \cite[Theorem 2]{BombieriMumford2}))
 
 There are integers $0 \leq a_i \leq m_i - 1$ and $\gamma_i$ such that
\begin{eqnarray*}
\omega_X & \cong& f^*(\omega_C \otimes \cal{L}^{-1}) \otimes  \OO_X(\sum_{i=1}^n a_iF_i) \\
m_i & = & p^{\gamma_i} \nu_i
\end{eqnarray*}
where $\deg(\omega_C \otimes \cal{L}^{-1}) = 2g(C) - 2 + \chi(X,\OO_X) + {\rm length}(\cal{T})$ and $\nu_i$ is the order of the normal bundle of $F_i$ in $X$.

\item (\emph{Ogg's formula} (see \cite{Ogg})) 

Let $\Delta_f$ be the discriminant of $f$ and $v_P(\Delta_f)$ the order of vanishing of $\Delta_f$ at $P \in C$. Then,
$$
c_2(X) = \sum_{P \in C} v_P(\Delta_f).
$$
Moreover, if $F_P$ denotes the fiber over $P$, then $v_P(\Delta_f) = e(F_P) + \delta_{F_P}$, where $e(F_P)$ is the topological Euler characteristic of $F_P$ and $\delta_{F_P}$ is the Swan conductor of $F_P$. If $m$ is the number of components of $F_P$, then
$$
e(F_P) = \begin{cases}
          0 & \text{ if } (F_P)_{red}\text{ is smooth }, \\
          m & \text{ if } F_P \text{ is multiplicative}, \\
          m+1  & \text{ if } F_P \text{ is additive}.
         \end{cases}
$$
and
$$
\delta_{F_P} = \begin{cases}
          0 & \text{ if } p \neq 2,3 \text{ or } F_P \text{ multiplicative}, \\
          0 & \text{ if } p = 3 \text{ and } F_P \text{ of type } {\rm III},{\rm III}^* \text{ or } {\rm I}_n^*, \\
          0 & \text{ if } p = 2 \text{ and } F_P \text{ of type } {\rm IV} \text{ or } {\rm IV}^*,  \\    
          \geq 2 & \text{ if } p = 2 \text{ and } F_P \text{ of type } {\rm II} \text{ or } {\rm I}_n^* \text{ with } n \neq 1, \\
          \geq 1 & \text{ else }.
         \end{cases}.
$$
For the list of Swan conductors, see e.g. \cite[p. 67]{schuettShioda}.

\item (\emph{Comparison of $f$ and $J(f)$ }(see \cite{LiuLorenziniRaynaud}))

Let $P \in C$ and $F_P$ resp. $F'_P$ be the fibers of $f$ resp. $J(f)$ over $P$. Then,
\begin{enumerate}[(i)]
 \item $f$ and $J(f)$ have the same $j$-map, 
 \item $F_P$ and $F'_P$ are of the same type,
 \item $v_P(\Delta_f) = v_P(\Delta_{J(f)})$,
 \item $\delta_{F_P} = \delta_{F'_P}$,
  \item $c_2(X) = c_2(J(X))$.
\end{enumerate}
\end {itemize}

Finally, we introduce the notion of \emph{$\Aut_X^0$-movable} fiber, which will play an important r\^ole throughout this article. The letters $a,m,\gamma,$ and $\nu$ will have the same meaning as the corresponding letters with indices in the canonical bundle formula recalled above.

\begin{Definition}\label{movabledefn}
A fiber $mF$ of an elliptic surface $f:X \to C$ is called \emph{$\Aut_X^0$-$n$-movable}, if $\Aut_{(nF,X)}^0 \neq \Aut_X^0$. An $\Aut_X^0$-$m$-movable fiber is simply called \emph{$\Aut_X^0$-movable}. We say that $mF$ is \emph{$n$-movable by vector fields}, if $H^0(X,T_X\langle nF \rangle) \neq H^0(X,T_X)$.
\end{Definition}

Clearly, an $\Aut_X^0$-$n$-movable fiber is $\Aut_X^0$-$k$-movable for all $k \leq n$, and if $mF$ is $\Aut_X^0$-$n$-movable by vector fields, it is $\Aut_X^0$-$n$-movable. The following lemma shows that $\Aut_X^0$-movable fibers satisfy very special properties.

\begin{Lemma} \label{movable}
Let $mF$ be an $\Aut_X^0$-$n$-movable fiber of $f$ with $n \geq 1$. Then, 
\begin{enumerate}[(i)]
\item $\nu = 1$, $m = p^{\gamma}$ with $\gamma \geq 0$ and either $a > n$ with $p \mid a$ or $a = 0$. In particular, if $mF$ is $\Aut_X^0$-movable, then $a = 0$.
\item $F$ is smooth, or $p \in \{2,3\}$ and $F$ is of type ${\rm II}$. In the latter case, $mF$ is not $1$-movable by vector fields.
\end{enumerate}
\end{Lemma}
\prf
By Proposition \ref{criteria} (a), we have $H^0(F,N_{F/X}) \neq 0$. Since $N_{F/X}$ has degree $0$ on every component of $F$, we deduce $N_{F/X} = \OO_{F}$ and therefore $\nu = 1$ and $m = p^\gamma$ for some $\gamma \geq 0$. 

Suppose that $a$ is prime to $p$. Then, there exists $l \geq 0$ such that $F$ appears as a reduced irreducible component of the scheme-theoretic base locus of $|l(K_X - f^*K_C)|$. Since $\Aut_X^0$ acts naturally on this base locus and preserves its connected components, we obtain $\Aut_{(F,X)}^0 = \Aut_X^0$. Hence, $p \mid a$. Moreover, if $a \neq 0$, then $aF$ is an irreducible component of the scheme-theoretic fixed locus of $|K_X - f^*K_C|$, so that $a > n$. Since $a$ is bounded above by the multiplicity of $F$, we deduce that $a = 0$ if $F$ is $\Aut_X^0$-movable.

Since $(-2)$-curves are infinitesimally rigid in $X$, Proposition \ref{criteria} (a) shows that $F$ is integral. Next, if $F$ is of type ${\rm I}_1$, then $\Def_{F}$ is prorepresentable and thus so is $\Def_{F/X}'$ by Lemma \ref{prorep}. Then, Proposition \ref{criteria} (c) shows that $mF$ is not $\Aut_X^0$-$1$-movable. 
If $F$ is of type ${\rm II}$, then $mF$ is not $\Aut_X^0$-$1$-movable by vector fields by Remark \ref{criteriaremark}. Moreover, if $p \neq 2,3$, then $\Def_{F}$ admits a TR-hull (see Remark \ref{cuspexample}), so that $F$ is not $\Aut_X^0$-$1$-movable by Lemma \ref{prorep} and Proposition \ref{criteria} (c).
\qed

\begin{Remark}
In Example \ref{movableexample}, we will give examples of elliptic surfaces over an affine curve with an $\Aut_X^0$-movable fiber of type ${\rm II}$ in characteristic $2$ and $3$.
\end{Remark}

\vspace{2mm}

\subsection{Automorphism schemes of genus $1$ curves}
In this section, we recall the structure of the automorphism scheme of a curve $C$ of genus $0$ or $1$ over an algebraically closed field $k$. This is well-known if $\Char(k) = p \neq 2,3$ and we refer the reader to \cite[Proposition 6]{BombieriMumford3} for proofs in the case of the cuspidal cubic if $p = 2,3$.

\begin{Lemma} \label{autoschemes}
Let $C$ be a reduced, irreducible curve of arithmetic genus $0$ or $1$ over $k$. Then, the following hold:
\begin{enumerate}[(i)]
\item If $C \cong \PP^1$, then $\Aut_C \cong \rm{PGL}_2$.
\item If $C$ is an elliptic curve, then $\Aut_C^0 \cong C$.
\item If $C$ is a nodal cubic curve, then $\Aut_C^0 \cong \GG_m$.
\item If $C$ is a cuspidal cubic curve, then $\Aut_C^0 \cong (\GG_a \rtimes A_p) \rtimes \GG_m$, where
$$
A_p = \begin{cases}
\{1\} & \text{ if } p \neq 2,3, \\
 \alpha_3 & \text{ if } p = 3, \\
 (\alpha_2 \times \alpha_2) \cdot \alpha_2 & \text{ if } p=2.
\end{cases}
$$ 

\end{enumerate}
\end{Lemma}

We will need further information on the fixed loci of some finite subgroup schemes of $\Aut_C^0$ in the above cases.

\begin{Lemma} \label{fixedloci}
Let $C$ be a reduced, irreducible curve of arithmetic genus $0$ or $1$ over $k$ and let $G \subseteq \Aut_C$ be a non-trivial connected subgroup scheme.
\begin{enumerate}[(i)]
\item If $C = \PP^1$, then $G$ has at most $2$ fixed points on $C$. Moreover, $G$ has precisely $2$ fixed points if and only if $G \subseteq \GG_m$.
\item If $C$ is an elliptic curve, then $G$ admits no fixed points on $C$.
\item If $C$ is a nodal cubic curve, then $G$ has exactly $2$ fixed points on $C$ and one of them is the node of $C$.
\item If $C$ is a cuspidal cubic curve and $G \cong \mu_{p^n}$, then one of the following holds
\begin{enumerate}[(1)]
\item $G$ has exactly $2$ fixed points on $C$ and one of them is the cusp of $C$,
\item $p^n = 2$, $G$ has exactly $3$ fixed points on $C$ and one of them is the cusp of $C$.
\item $p^n = 2$, $G$ has exactly $4$ fixed points on $C$ and all of them are smooth points of $C$,
\item $p^n = 3$, $G$ has exactly $3$ fixed points on $C$ and all of them are smooth points of $C$,
\item $p^n = 4$, $G$ has exactly $2$ fixed points on $C$ and both of them are smooth points of $C$. In this case, the induced $\mu_2$-action is as in case (2).
\end{enumerate}
\end{enumerate}
\end{Lemma}

\prf
Claims (ii) and (iii) and the first part of Claim (i) are well-known. Let us prove the second part of Claim (i). If $G$ fixes two points on $\PP^1$, then $G \subseteq \GG_m$. Conversely, if $G \subseteq \GG_m$, then we can conjugate $G \subseteq {\rm PGL}_2$ so that it lies in the diagonal torus. Then, $G$ fixes $0$ and $\infty$ on $\PP^1$.

To prove Claim (iv), we recall that by \cite[Proposition 6]{BombieriMumford3}, one can identify the smooth locus of $C$ with $\mathbb{A}^1 = \Spec k[t]$ such that automorphisms of $\mathbb{A}^1$ induced by automorphisms of $C$ are of the following form:
\begin{eqnarray*}
t \mapsto& at + b,  &a \in \GG_m, b \in \GG_a \text{ if } p \neq 2,3, \\
t \mapsto& at + b + ct^3, & a \in \GG_m, b \in \GG_a, c^3 = 0 \text{ if } p = 3, \\
t \mapsto& at + b + ct^2 + dt^4,  & a \in \GG_m, b \in \GG_a, c^4 = d^2 = 0 \text{ if } p = 2.
\end{eqnarray*}
Moreover, we refer the reader to \cite[p. 212]{BombieriMumford3} for the calculation of the stabilizer of the cusp of $C$, which is given by all substitutions if $p \neq 2,3$, by the substitutions with $c = 0$ if $p = 3$, and by the substitutions with $c^2 = d = 0$ if $p = 2$.

If $p \neq 2,3$, then $G$ is conjugate to the $\mu_{p^n}$ of maps $t \mapsto at$, $a \in \mu_{p^n}$. Its fixed points are $t = 0$ and the cusp of $C$.

If $p = 3$, then we can conjugate $G$ such that either $G$ acts as above or as $t \mapsto at + (1-a)t^3$ with $a^3 = 1$. In the latter case, the fixed points are given by $t^3 = t$. This $\mu_3$-action does not fix the cusp of $C$.

If $p = 2$, we can conjugate $G$ such that it acts in one of the following ways with $\lambda,\mu \in k$:
\begin{eqnarray*}
t  \mapsto & at, & a \in \mu_{2^n}, \\
t  \mapsto & at + \lambda(1+a)t^2 + \mu(1+a)t^4, & a \in \mu_2, \\
t  \mapsto & at + (a+a^2)t^2 + (1+a^2)t^4,&  a \in \mu_4.
\end{eqnarray*}
In the first case, $G$ fixes $t = 0$ and the cusp of $C$. 
In the second case, $G$ fixes the points where $\mu t^4 + \lambda t^2 + t = 0$. If $\mu \neq 0$, this shows that $G$ fixes $4$ smooth points on $C$ and does not fix the cusp, whereas if $\mu = 0$, the action of $G$ has $2$ smooth fixed points on $C$ and fixes the cusp.
In the third case, $G$ fixes $t \in \{0,1\}$ and does not fix the cusp of $C$. Moreover, $G[F] = \mu_2$ acts as in the second case with $\mu = 0$. This proves Claim (iv).
\qed
\vspace{3mm}

Using Lemma \ref{fixedloci}, we can determine how $\mu_{p^n}$-actions on an elliptic surface $f: X \to C$ can restrict to reducible fibers of $f$ and determine the possible fixed loci.

\begin{Proposition}\label{fixedlocireducible}
Let $f:X \to C$ be an elliptic surface with $\mu_{p^n} \subseteq \Aut_X^0$. Let $mF$ be a singular fiber of $f$. Then, the following hold:
\begin{enumerate}[(i)]
\item If $F$ is not of type ${\rm II}$ or ${\rm III}$, then $e(F) = e(F^{\mu_{p^n}})$.
\item If $F$ is of type ${\rm III}$, then $e(F) = e(F^{\mu_{p^n}}) = 3$, or $p^n = 2$ and $e(F^{\mu_2}) = 4$.
\item If $F$ is of type ${\rm II}$, then $F$ is preserved by the $\mu_{p^n}$-action, and
\begin{enumerate}[(1)]
\item $e(F) = e(F^{\mu_{p^n}}) = 2$, or 
\item $p^n = 3$ and $e(F^{\mu_3}) = 3$, or 
\item $p^n = 2$ and $e(F^{\mu_2}) \in \{3,4\}$.
\end{enumerate}
\end{enumerate}
\end{Proposition}

\prf
If $F$ is not of type ${\rm II}$, then $mF$ is not $\Aut_X^0$-$1$-movable by Lemma \ref{movable}, so the $\mu_{p^n}$-action on $X$ restricts to a $\mu_{p^n}$-action on $F$. 

In the first case, all intersections of components of $F$ are transversal. Since $\mu_{p^n}$ preserves all components, it fixes all their intersections. Now, the statement can be checked case by case and the proof is the same as in \cite[Lemma 2]{DolgachevNumerical}.

In the second case, the intersection of the two components $E_1,E_2$ of $F$ is not transversal. If $\mu_{p^n}$ fixes $(E_1 \cap E_2)_{red}$, then $e(F^{\mu_{p^n}}) = 3$ by Lemma \ref{fixedloci} (i) and, since $e(F) = 3$, this gives the desired equality of Euler characteristics. If $\mu_{p^n}$ does not fix $(E_1 \cap E_2)_{red}$, then it fixes two points on each of the $E_i$ by Lemma \ref{fixedloci} (i). Consider the contraction $\pi: X \to X'$ of $E_1$. Then, $\pi(E_2)$ is a cuspidal curve on $X'$, the $\mu_{p^n}$-action on $X'$ induced via Theorem \ref{Blanchard} has three fixed points on $\pi(E_2)$ and one of them is the cusp of $\pi(E_2)$ by Proposition \ref{blowup}. Hence, $p^n = 2$ by Lemma \ref{fixedloci}.

If $F$ is of type ${\rm II}$, the only statement that is not already included in Lemma \ref{fixedloci} is the fact that $F$ is preserved by $\mu_{p^n}$. To prove this, note that $F$ is not $1$-movable by vector fields by Lemma \ref{movable} and hence $\mu_{p^n} \cap \Aut_{(F,X)}^0$ is non-trivial. In particular, the induced $\mu_p$-action preserves $F$. By Lemma \ref{fixedloci}, this $\mu_p$-action has an isolated fixed point $Q$ on $F$.  Since $X^{\mu_p}$ is smooth at $Q$, the point $Q$ is also a fixed point of the $\mu_{p^n}$-action by Proposition \ref{isolatedefixedpoints}. Therefore, by Proposition \ref{blowup}, the $\mu_{p^n}$-action lifts to the blowup $\widetilde{X}$ of $X$ at $Q$. Since the strict transform $\widetilde{F}$ of $F$ in $\widetilde{X}$ is a negative curve, we have $\Aut_{(\widetilde{F},\widetilde{X})}^0 = \Aut_{\widetilde{X}}^0$ by Proposition \ref{criteria} and therefore the $\mu_{p^n}$-action on $X$ preserves $F$ by Lemma \ref{image}.
\qed

\begin{Remark}
It was claimed in \cite[Lemma 3]{RudakovShafarevich} that the exceptional case in Proposition \ref{fixedlocireducible} (ii) does not occur. We will give a counterexample to this statement in Example \ref{caseb2}. The proof of \cite[Lemma 3]{RudakovShafarevich} seems to be correct up until the last sentence, where it is claimed that the configuration described in \cite[p. 1224]{RudakovShafarevich} is not of Kodaira type. In fact, the configuration described there is of type ${\rm I}_1^*$.
\end{Remark}

\begin{Remark}  \label{fixedlocireducibleremark}
The proof of Proposition \ref{fixedlocireducible} (iii) shows more generally that an irreducible fiber $F$ of an elliptic fibration $f: X \to C$ is preserved by a group scheme action as soon as the action has a fixed point on $F$. 
\end{Remark}

\begin{Example} \label{movableexample}
The following examples show that, at least locally, there may be group scheme actions on elliptic surfaces that actually move fibers of type ${\rm II}$ in characteristic $2$ and $3$: Let $p = 2$ and let $X \subseteq  \PP^2_{k[t]}$ be the smooth surface defined by
$$
y^2z + t^4yz^2 = x^3 + tz^3.
$$
The generic fiber of $X \to \Spec k[t]$ is an elliptic curve and the fiber $F$ at $t = 0$ is a cuspidal cubic.
There is an $\alpha_4$-action on $X$ defined by
\begin{eqnarray*}
(x,z) &\mapsto& (x,z)  \\
y &\mapsto& y + az\\
t &\mapsto& t + a^2 + at^4
\end{eqnarray*}
where $a^4 = 0$. Note that the induced $\alpha_2$-action preserves $F$, but the $\alpha_4$-action itself does not.
Moreover, the $\alpha_4$-action has no fixed point on $f$, since the induced $\alpha_4$-action on $\Spec k[t]$ has no fixed point.
 A similar example of an $\alpha_9$-action in characteristic $3$ exists on the surface defined by
$$y^2z = x^3 + t^9xz^2 + tz^3.$$
\end{Example}

\section{Vertical component of $\Aut_X^0$} \label{vertical}
In this section, $f: X \to C$ is an elliptic surface over an algebraically closed field $k$. The purpose of this section is to study the vertical component of $\Aut_X^0$.
Recall that because of Blanchard's Lemma (see Theorem \ref{Blanchard}) there is a natural map $f_*: \Aut_X^0 \to \Aut_C^0$.

\begin{Definition}
Let $f: X \to C$ be an elliptic surface. The \emph{vertical component} of $\Aut_X^0$ is defined as $\Ker(f_*)$, where $f_*: \Aut_X^0 \to \Aut_C^0$ is the natural map.
\end{Definition}

After recalling the notion of Weil restrictions of group schemes along the field extension $h: \Spec k(C) \to \Spec k$, we will first study the action of $\Ker(f_*)$ on the generic fiber $F_\eta$ of $f$ and then determine obstructions to extending such actions to the surface $X$.

\subsection{Automorphisms of the generic fiber}

Recall the following results on Weil restrictions from \cite[Section 7.6]{NeronModels}.

\begin{Definition}
The \emph{Weil restriction} of a scheme $G$ over $k(C)$ along $h: \Spec k(C) \to \Spec k$ is defined as the presheaf
\begin{eqnarray*}
h_*G: & (Sch/k)^{op} &\to (Sets) \\
& T &\mapsto G(T \times_{\Spec k} \Spec k(C)).
\end{eqnarray*}
\end{Definition}

\begin{Lemma}\label{Weilrestrictionlemma}
There is a bijection of sets of homomorphisms of presheaves
$$
\Hom_{\Spec k}(T,h_*G) \to \Hom_{\Spec k(C)}(T \times_{\Spec k} \Spec k(C), G)
$$
which is functorial in the $k$-scheme $T$ and the $k(C)$-scheme $G$.
\end{Lemma}

\begin{Lemma} \label{unexpectedkernel}
Let $G$ be a group scheme over $k(C)$, let $G'$ be a group scheme over $k$ and let $g: G'_{k(C)} \to G$ be a morphism of group schemes such that the induced map $h_*g: G' \to h_* G$ is a monomorphism of presheaves. Then, the only subscheme of $\Ker(g)$ that can be defined over $k$ is the trivial subgroup scheme. In particular, $\Ker(g)_{red}$ is trivial.
\end{Lemma}

\prf
Let $H \subseteq \Ker(g)$ be a subscheme which can be defined over $k$. Then, by definition, there exists a scheme $H'$ over $k$ and a morphism $H' \to G'$ whose base change along $h$ agrees with $H \to G'_{k(C)}$. Since the induced map $H \to G$ is constant, it follows from the adjunction in Lemma \ref{Weilrestrictionlemma} that the map $H' \to h_*G$ is constant. But $h_*g$ is a monomorphism, hence $H'$ is trivial and thus so is $H$. In particular, $\Ker(g)$ contains only one point and hence $\Ker(g)_{red}$ is trivial.
\qed
\vspace{3mm}

We will now apply the Weil restriction to automorphisms of elliptic surfaces. For every $k$-scheme $T$, we have a natural injective map
$$
\Ker(f_*)(T) \to \Aut_{F_\eta/k(C)}(T \times_{\Spec k} \Spec k(C)),
$$
where $F_\eta$ is the generic fiber of $f$. This defines a monomorphism of presheaves of groups $\Ker(f_*) \to h_*\Aut_{F_\eta/k(C)}$ and hence we obtain a morphism of group schemes $\varphi: \Ker(f_*)_{k(C)} \to \Aut_{F_\eta/k(C)}$ from Lemma \ref{Weilrestrictionlemma}. The connected component of the identity of the latter group scheme is isomorphic to the generic fiber $J_\eta$ of the Jacobian $J(f)$ of $f$ and we denote the induced map of identity components by $\varphi^0: \Ker(f_*)_{k(C)}^0 \to J_\eta$. The following lemma shows that $\varphi^0$ is injective as long as $F_\eta$ is ordinary.

\begin{Lemma}\label{connsub}
Let $f: X \to C$ be an elliptic surface. Let $G = \Ker(f_*)^0$ and let $\varphi^0$ as above.
\begin{enumerate}[(i)]
\item The group scheme $G$ is commutative and $\dim(G) \leq 1$,
\item If $\dim(G) = 1$, then the Jacobian $J(f): J(X) \to C$ of $f$ is a trivial fibration.
\item If $J_\eta$ is ordinary, then $\varphi^0$ is injective. In this case, either $\dim(G) = 1$ and $G_{k(C)} \cong J_\eta$, or $\dim(G) = 0$ and $G \cong \mu_{p^n}$ for some $n \geq 0$.
\item If $J_\eta$ is supersingular, then we have $G[F] \cong \alpha_p^r$ for some $r \geq 0$. If $r = 1$, then either $\dim(G) = 1$ and $G_{k(C)} \cong J_\eta$, or $\dim(G) = 0$ and $G \cong M_n$ for some $n \geq 0$.\end{enumerate}
\end{Lemma}

\prf
First, note that the action of $G$ on $X$ factors through $G^{ab}$ on a dense open subset of $X$, because the identity component of the automorphism scheme of a smooth curve of genus one is commutative. Therefore, the action of $G$ on all of $X$ factors through $G^{ab}$ and thus $G$ is commutative, proving the first part of (i).

To prove Claim (i) and (ii), let $G_r := G_{red}$ be the reduction of $G$. Since $k$ is perfect, Lemma \ref{generalgroupscheme} (iv) shows that this is a closed and smooth subgroup scheme of $G$. Assume that $\dim(G) \geq 1$. Then, $G_r$ is non-trivial. Consider the morphism $\varphi_r: (G_r)_{k(C)} \to J_\eta$ obtained by restricting $\varphi^0$ to $(G_r)_{k(C)}$. By Lemma \ref{unexpectedkernel}, the group scheme $\Ker(\varphi_r)$ is zero-dimensional and connected, hence $\varphi_r$ is a purely inseparable isogeny of elliptic curves over $k(C)$. But all finite connected subgroup schemes of $(G_r)_{k(C)}$ are of the form $(G_r)_{k(C)}[F^n]$ for some $n \geq 0$. Since $(G_r)_{k(C)}[F^n] = (G_r[F^n])_{k(C)}$,  these subschemes can be defined over $k$. Hence, by Lemma \ref{unexpectedkernel}, the map $\varphi_r$ is an isomorphism. Then, $G \times C \to C$ is a minimal proper regular model for $J_\eta$ over $C$ and hence coincides with $J(f)$ by Lemma \ref{uniquemodel}. In particular, $J(f)$ is a trivial fibration. This yields Claim (i) and (ii).

To prove Claim (iii) and (iv), we use Lemma \ref{finitecommgroupscheme} (ii) to write $G[F^n] \cong G[F^n]_{lr} \times G[F^n]_{ll}$ and consider the action of $G[F^n]$ on $X$.

For Claim (iii), assume that $J_\eta$ is ordinary. Then, almost all fibers of $f$ are ordinary. By Lemma \ref{homsandexts} (i), every action of $G[F^n]_{ll}$ on an ordinary elliptic curve is trivial, hence $G[F^n]$ is of type $lr$ and thus isomorphic to $\prod_{i=1}^m \mu_{p^{n_i}}$ for some $m,n_i \geq 0$ by Lemma \ref{finitecommgroupscheme} (ii) (3). Therefore, we have $(G[F^n])_{k(C)} \cong \prod_{i=1}^m (\mu_{p^{n_i}})_{k(C)}$. Subgroup schemes of this group scheme correspond to quotients of its reduced Cartier dual, hence all of them are defined over $k$. Thus, by Lemma \ref{unexpectedkernel} the intersection $\Ker(\varphi^0) \cap G[F^n]$ is trivial and we have $(G[F^n])_{k(C)} \cong J_\eta[F^n] \cong (\mu_{p^{n}})_{k(C)}$ and thus $G[F^n] \cong \mu_{p^n}$. Since $\Ker(\varphi^0)$ is finite and connected, we have $\Ker(\varphi^0) \subseteq G[F^n]$ for $n \gg 0$, so we can in fact deduce that $\varphi^0$ is injective. If $\dim(G) = 0$, then $G \cong G[F^n]$ for $n \gg 0$ and if $\dim(G) = 1$, then $G = G_{red}$ and, as in the second paragraph of the proof, $\varphi^0$ induces an isomorphism $G_{k(C)} \to J_\eta$. This yields Claim (iii).

As for Claim (iv), we assume that $J_\eta$ is supersingular. Then, the group $G[F]$ is of type $ll$ by Lemma \ref{homsandexts} (i). Moreover, the action of $G[F]$ on a general fiber $E$ of $f$ factors through $E[F]$. Since Verschiebung is trivial on $E[F]$, the action of $G[F]$ on $X$ factors through $G[F]/(VG[F])$, hence $VG[F] = 0$. Therefore, by Lemma \ref{finitecommgroupscheme} (ii), we have $G[F] \cong \alpha_p^r$ for some $r \geq 0$. Now, if $r = 1$, then $\alpha_p$ is the unique simple closed subgroup scheme of $G[F^n]$ for every $n \geq 1$. Therefore, if the morphism $G[F^n] \to E[F^n]$ is not injective for a general fiber $E$, then $G[F] = \alpha_p$ is in its kernel and therefore $(G[F])_{k(C)} \subseteq \Ker(\varphi)$, which is impossible by Lemma \ref{unexpectedkernel}. Hence, $G[F^n]$ is isomorphic to its image in $E[F^n] \cong M_n$. This yields Claim (iv).
\qed
\vspace{2mm}

\begin{Remark}
In Example \ref{importantexample}, we will show that the integer $r$ appearing in Lemma \ref{connsub} (iv) can be arbitrarily large in every positive characteristic.
\end{Remark}

If $J_\eta$ is ordinary, the existence of a subscheme isomorphic to $(\mu_{p^n})_{k(C)}$ with $n \geq 1$ in $J_\eta$ already gives strong restrictions on the geometry of $J(f)$. This is closely related to the Igusa moduli problem, which is defined as follows.

\begin{Definition} \label{Igusadefn}
The ordinary part ${\rm Ig}(p^n)^{\rm ord}$ of the \emph{Igusa stack} is the stack over the category of $k$-schemes whose objects over a $k$-scheme $T$ are families $E \to T$ of ordinary elliptic curves over $T$ together with a generator of $E^{(p^n)}[V^n] := {\rm Ker}(V^n)(E^{(p^n)} \to E)$, where $V: E^{(p)} \to E$ denotes Verschiebung on $E$.
\end{Definition}

This moduli problem has been first studied by Igusa in \cite{Igusa}. If $p^n \geq 3$, then ${\rm Ig}(p^n)^{\rm ord}$ is representable by a smooth curve defined over $\FF_p$ (see \cite[Corollary 12.6.3]{KatzMazur}). We denote its smooth projective compactification by ${\rm Ig}(p^n)$.
Now, the following lemma is a straightforward consequence of the definition of ${\rm Ig}(p^n)$.

\begin{Lemma} \label{Igusa}
Let $J(f): J(X) \to C$ be a Jacobian elliptic surface with generic fiber $J_\eta$. Assume that $p^n > 2$. Every inclusion $(\mu_{p^n})_{k(C)} \hookrightarrow J_\eta$ yields a morphism $u: \Spec k(C) \to {\rm Ig}(p^n)$ such that $J_\eta$ is the pullback of the universal elliptic curve over ${\rm Ig}(p^n)$ along $u$.
\end{Lemma}

\prf
Since $(\mu_{p^n})_{k(C)}$ has length $p^n$, the inclusion $(\mu_{p^n})_{k(C)} \hookrightarrow J_\eta$ factors through an isomorphism $(\mu_{p^n})_{k(C)} \cong J_\eta[F^n]$. Taking Cartier duals and using $(\mu_{p^n})_{k(C)}^{\vee} = \mathbb{Z}/p^n\mathbb{Z}$, we obtain an isomorphism
$
J_{\eta}^{(p^n)}[V^n] \cong (J_\eta[F^n])^{\vee} \cong \mathbb{Z}/p^n\mathbb{Z}.
$
The preimage of $1$ under this isomorphism yields a generator of $J_{\eta}^{(p^n)}[V^n]$. Hence, by Definition \ref{Igusadefn}, we get a morphism $u: \Spec k(C) \to {\rm Ig}(p^n)$ inducing $J_\eta$ via pullback of the universal elliptic curve over ${\rm Ig}(p^n)$.
\qed
\vspace{3mm}

In \cite{LiedtkeSchroer}, Liedtke and Schr\"oer studied the singular fibers of the universal elliptic surfaces over the Igusa curves ${\rm Ig}(p)$. Using their results and Lemma \ref{Igusa}, we obtain the following proposition.

\begin{Proposition}\label{singularfibers}
Assume $p > 3$. Let $f: X \to C$ be an elliptic surface with $(\mu_p)_{k(C)} \subseteq \Aut_{F_\eta/k(C)}^0$. Let $mF$ be an additive fiber of $f$. Then, the following hold:
\begin{enumerate}[(i)]
\item $F$ is not of type ${\rm I}_n^*$ with $n \geq 1$.
\item If $p \equiv 1$ mod $12$, then $F$ is of type ${\rm I}_0^*$.
\item If $p \equiv 7$ mod $12$, then $F$ is of type ${\rm III},{\rm III}^*$ or ${\rm I}_0^*$.
\item If $p \equiv 5$ mod $12$, then $F$ is of type ${\rm II},{\rm IV},{\rm IV}^*,{\rm II}^*$ or ${\rm I}_0^*$.
\end{enumerate}
\end{Proposition}
\prf
We have $(\mu_p)_{k(C)} \subseteq \Aut_{F_\eta/k(C)}^0 = J_\eta$, where $J_\eta$ is the generic fiber of the Jacobian of $f$, identified with $\Aut_{F_\eta/k(C)}^0$ via its natural action on $F_{\eta}$. Therefore, the curve $J_\eta$ is a pullback of the universal elliptic curve over ${\rm Ig}(p^n)$ along a morphism $\Spec k(C) \to {\rm Ig}(p^n)$  by Lemma \ref{Igusa}. Since $f$ and $J(f)$ have the same types of singular fibers, the claim now follows by comparing the reduction types of the universal elliptic curve given in \cite[Theorem 10.1 and Theorem 10.3]{LiedtkeSchroer} with the tables in \cite[Section 5.2.]{schuettShioda}.
\qed

\begin{Remark}
There is no analogue of Lemma \ref{Igusa} if $J_\eta$ is supersingular. For $\alpha_p$ and $M_2$, this follows immediately from \cite[Theorem 6.1]{Liedtkeptors}, and for $M_n$ with $n \geq 3$ one can simply iterate the argument given in the proof there.
\end{Remark}

The genus $g({\rm Ig}(p^n))$ of the Igusa curve has been computed by Igusa in \cite{Igusa}. We have
$$
g({\rm Ig}(p^n))  = \frac{1}{48}(p-1)(p^{2n-1}-12p^{n-1} + 1) + 1 - \frac{h_p}{2},
$$
where 
$$
h_p = \text{number of supersingular $j$-invariants in } k +
\begin{cases}
0  & \text{ if } p \neq 2,3 \\
\frac{1}{3} &\text{ if } p = 3 \\
\frac{3}{8} & \text{ if } p = 2
\end{cases}.
$$  In particular,
$$
g({\rm Ig}(p^n)) = 
\begin{cases}
 0 &\text{ if } p^n \leq 12 \\
 1 &\text{ if } p^n \in \{13,16\}  \\
 \geq 2 &\text{ else}
\end{cases}
$$

\begin{Proposition}\label{Igusainequality}
Let $f: X \to C$ be an elliptic surface with $c_2(X) \neq 0$. If $\mu_{p^n} \subseteq \Ker(f_*)$, then
$$
g(C) \geq \frac{1}{48}(p-1)(p^{2n-1}-12p^{n-1} + 1) + 1 - \frac{h_p}{2}
$$
\end{Proposition}

\prf
We can assume that $p^n > 2$, for otherwise the right hand side of the inequality is negative. Since the genus of smooth curves does not go down under taking finite covers, it suffices to show that $C$ admits a dominant rational map to ${\rm Ig}(p^n)$. 

By Lemma \ref{connsub} (iv), the map $\varphi: (\mu_{p^n})_{k(C)} \to J_\eta$ induced by the inclusion $\mu_{p^n} \subseteq \Ker(f_*)$ is injective. Hence, Lemma \ref{Igusa} shows that there is a morphism $\Spec k(C) \to {\rm Ig}(p^n)$. Seeking a contradiction, we assume that this map is constant. Then, we have $J_\eta = E \times_{\Spec k} \Spec k(C)$ for some ordinary elliptic curve $E$ over $k$. Thus, $E \times C$ is a minimal proper regular model of $J_\eta$ over $C$ and therefore it coincides with $J(X)$ by Lemma \ref{uniquemodel}. Then, $0 = c_2(E \times C) = c_2(J(X)) = c_2(X)$, contradicting our assumption. Hence, $\Spec k(C) \to {\rm Ig}(p^n)$ is dominant, which is what we had to prove.
\qed
\vspace{5mm}

\subsection{Extending the action to $X$}
In the previous subsection, we have seen how the existence of a $\mu_{p^n}$-action on the generic fiber $F_\eta$ of an elliptic surface $f:X \to C$ gives restrictions on $f$. In this section, we gather several criteria for a connected group scheme action on $F_\eta$ to extend to an action on $X$. Using these criteria, we give geometric restrictions that have to be satisfied by elliptic surfaces with non-trivial $\Ker(f_*)$. 
For the following definition, note that the $k(C)$-scheme $F_{\eta}$ is a $C$-scheme via $F_{\eta} \to \Spec k(C) \to C$ and for every $k$-group scheme $G$ we can consider $G \times F_{\eta}$ as a $C$-scheme via the second factor.

\begin{Definition} \label{extendabilitydefinition}
Let $f: X \to C$ be an elliptic surface with generic fiber $F_\eta$. Let $G$ be a $k$-group scheme and let $\rho_\eta: G \times F_\eta \to F_\eta$ be an action of $G$ on $F_\eta$ over $k$ such that $\rho_\eta$ is a morphism of $C$-schemes. We say that $\rho_\eta$ \emph{extends over} $p \in C$ if there is a commutative diagram
$$
\xymatrix{
G \times F_\eta \ar[r]^{\rho_{\eta}} \ar[d] & F_\eta \ar[d] \\
G \times X_p \ar[r]^{\rho_p}  & X_p,
}
$$
where $X_p := (X \times \Spec \OO_{C,p})$ and the vertical arrows are induced by the inclusion $\OO_{C,p} \subseteq k(C)$.
We say that $\rho_\eta$ \emph{extends to} $X$ if there is a similar diagram with a morphism $\rho: G \times X \to X$ in the second row.
\end{Definition}

Equivalently, via the natural isomorphism $G \times F_{\eta} \cong G_{k(C)} \times_{\Spec k(C)} F_{\eta}$, we can think of $\rho_{\eta}$ as an action of $G_{k(C)}$ on $F_{\eta}$ over $k(C)$ and ask whether \mbox{it extends to an action of $G_{\OO_{C,p}}$ on $X_p$ over $\OO_{C,p}$.}

\begin{Remark}
Note that $\rho$ and $\rho_p$ are automatically actions of $G$, as associativity can be checked on the schematically dense subscheme $G \times G \times F_\eta$ of $G \times G \times X$.
\end{Remark}

Recall from the previous subsection that an action $\rho_\eta: G \times F_\eta \to F_\eta$ as above gives rise to a translation action $\rho'_\eta: G \times J_\eta \to J_\eta$, where $J_\eta$ is the generic fiber of the Jacobian $J(f)$ of $f$. In the following proposition, we relate extendability of $\rho_\eta$ to extendability of $\rho'_\eta$.

\begin{Proposition} \label{extending}
Let $f: X \to C$ be an elliptic surface with generic fiber $F_\eta$ and let $\rho_\eta: G \times F_\eta \to F_\eta$ be an action of a connected group scheme $G$ such that $\rho_\eta$ is a morphism of $C$-schemes. Then, the following hold:
\begin{enumerate}[(i)]
\item The action $\rho_\eta$ extends to $X$ if and only if it extends over every $p \in C$.
\item If $p \in C$ is a point such that the fiber $F_p$ of $f$ over $p$ is simple, then $\rho_\eta$ extends over $p$ if and only if the corresponding action $\rho'_\eta$ on the generic fiber $J_\eta$ of the Jacobian $J(f)$ of $f$ extends over $p$.
\end{enumerate}
\end{Proposition}

\prf
The action $\rho_\eta$ gives rise to a rational map $C \dashrightarrow {\rm Hom}_C(G \times F_\eta, F_\eta)$. Since the latter scheme is separated and $C$ is a smooth curve, this rational map extends to a morphism if and only if it extends over every closed point of $C$. This shows Claim (i).

Next, let $p \in C$ be a point such that $F_p$ is simple and let $A := \Spec \OO_{C,p}$. Since the smooth locus of $X_p \to A$ is a torsor under its Jacobian and $\rho_\eta$ is induced by restricting the action of $J_\eta$ on $F_\eta$ to $G$, there is an \'etale cover $B \to A$, which we may assume to be Galois with covering group $H$, and a $G$-equivariant isomorphism of $B$-schemes
$$
\varphi: J(X)_p \times_A B \cong  X_p \times_A B.
$$
Both sides of the isomorphism are equipped with the natural action of $H$ on the second factor and the $G$-action on $J(X)_p \times_A B$ (resp. $X_p \times_A B$) descends to $J(X)_p$ (resp. $X_p$) if and only if it is normalized by this action of $H$. Thus, we obtain two actions of $H$ on both sides of the above isomorphism and one can check that these two actions differ by translation by a $B$-valued section of $J(X)_p \times_A B$ (see \cite[p.1233]{RudakovShafarevich}). By construction, $\rho_\eta$ and $\rho'_\eta$ commute with translations, so if one of them, say $\rho_\eta$, extends over $p$, then the induced $G$-action on $X_p \times_A B$ is normalized by both $H$-actions, hence this $G$-action also descends to $J(X)_p$. The induced action agrees with $\rho'_\eta$ on $J_\eta$, hence it extends $\rho'_\eta$ over $p$. By the same argument, $\rho_\eta$ extends over $p$ if $\rho_\eta'$ does. This proves Claim (ii).
\qed
\vspace{5mm}

Hence, if $f: X \to C $ admits no multiple fibers, then $G$-actions on $X$ which are trivial on $C$ correspond naturally to $G$-actions on $J(X)$ which are trivial on $C$, that is, we have the following corollary.
\begin{Corollary}\label{ExtendabilityCorollary}
Let $f: X \to C$ be an elliptic surface without multiple fibers. Then, $\Ker(f_*)^0 \cong \Ker(J(f)_*)^0$.
\end{Corollary}

The situation becomes more complicated if $f$ admits multiple fibers. Nevertheless, it turns out that an elliptic fibration $f$ with non-trivial $\Ker(f_*)^0$ must satisfy severe geometric constraints.

\begin{Theorem}\label{Main}
Let $f:X \to C$ be an elliptic surface with $\Ker(f_*)^0$ non-trivial. Then, the following hold:
\begin{enumerate}[(i)]
\item Every separable multisection $\Sigma$ of $f$ satisfies $\Sigma^2 \geq 0$.
\item Either $\chi(X,\OO_X) = 0$ or $f$ admits a multiple fiber.
\item One of the following two cases holds:
\begin{enumerate}[(1)]
\item $\alpha_p^r \subseteq \Ker(f_*)$ for some $r \geq 1$, every separable multisection $\Sigma$ of $f$ satisfies the inequality $h^0(\Sigma,N_{\Sigma/X}) \geq r$ and $f$ is isotrivial with supersingular generic fiber, or
\item $\mu_{p^n} \subseteq \Ker(f_*)$ for some $n \geq 1$ and all additive or supersingular fibers of $f$ are multiple fibers with multiplicity divisible by $p^n$.
\end{enumerate}
\end{enumerate}
\end{Theorem}

\prf
Let $\Sigma \subseteq X$ be an irreducible curve such that $f|_\Sigma: \Sigma \to C$ is finite and separable. Then, the curve $\Sigma$ intersects a general fiber of $f$ transversally, say in $n$ points. Since $\Ker(f_*)^0$ acts without fixed point on a general fiber of $f$ by Lemma \ref{fixedloci} (ii), this implies that $\Sigma$ is not preserved by $\Ker(f_*)^0$. Thus, by Proposition \ref{criteria}, we have $\Sigma^2 \geq 0$. This is Claim (i).

To prove Claim (ii), we need to show that if $f$ admits no multiple fibers, then $\chi(X,\OO_X) = 0$. By Corollary \ref{ExtendabilityCorollary} and since $\chi(X,\OO_X) = \chi(J(X),\OO_{J(X)})$, we may assume that $f$ admits a section $\Sigma$. Applying adjunction and the canonical bundle formula, we obtain
$$
2g(\Sigma) - 2 = 2g(C) - 2  + \chi(X,\OO_{X}) + \Sigma^2.
$$
Since $\chi(X,\OO_{X})$ is always non-negative and $\Sigma^2$ is non-negative by Claim (i), we deduce from $g(\Sigma) = g(C)$ that $\chi(X,\OO_X) = \Sigma^2 = 0$, as claimed.

For the proof of Claim (iii), let us first assume that $\alpha_p^r \subseteq \Ker(f_*)$ for some $r \geq 1$. If $h^0(\Sigma,N_{\Sigma/X}) \leq r-1$, then $\Aut_{(\Sigma,X)}^0 \cap \Ker(f_*)$ is non-trivial, which is impossible by the same argument as in the first paragraph, since $\Sigma$ intersects a general fiber of $f$ transversally. By Lemma \ref{connsub}, the existence of $\alpha_p \subseteq \Ker(f_*)$ implies that the generic fiber of the Jacobian of $f$ is supersingular. Since a supersingular elliptic curve can be defined over a finite field, this implies that $J(f)$ is isotrivial and therefore the same holds for $f$.

If $\alpha_p \not \subseteq \Ker(f_*)$, then $\mu_{p^n} \subseteq \Ker(f_*)$ for some $n \geq 1$ by Lemma \ref{connsub}. Let $mF$ be an additive or supersingular fiber of $f$, where $m$ is the multiplicity of $F$ and let $P = f(F)$. To finish the proof, we have to show that $p^n \mid m$.

First, assume that $\mu_{p^n}$ preserves $F$. Since $F$ is additive or smooth, there is a reduced component $F_1$ of $F$ which meets at most one other component of $F$. Then, Lemma \ref{autoschemes} shows that $\mu_{p^n}$ has a fixed point $Q$ on $F_1$ that does not lie on any other component of $F$.
By Section \ref{Linearization}, we can linearize the $\mu_{p^n}$-action in a formal neighborhood of $Q$, i.e. there is a $\mu_{p^n}$-equivariant isomorphism
$$
\widehat{\OO}_{X,Q} \cong k[[x,y]]
$$
such that $\mu_{p^n}$ acts via $x \mapsto \lambda x, y \mapsto \lambda^a y$ for some $0 \leq a \leq p^n-1$. 
Since $F_1$ is preserved by the $\mu_{p^n}$-action and the fixed locus of the $\mu_{p^n}$-action is contained in fibers of $f$, we can assume without loss of generality that $F_1$ is defined by $x = 0$. Let $t$ be a parameter on $C$ at $P$. The morphism $\varphi^{\#}: \widehat{\OO}_{C,P} \cong k[[t]] \to k[[x,y]]$ is then given by $\varphi^{\#}(t) = ux^m$, where $u \in k[[x,y]]$ is a unit and $m$ is the multiplicity of $F_1$.
Now, since the $\mu_{p^n}$-action on $C$ is trivial, we must have $ux^m \in k[[x,y]]^{\mu_{p^n}}$. In particular, the leading monomial of $ux^m$, which is of the form $c x^m$ for some $c \in k^{\times}$, has to be $\mu_{p^n}$-invariant. Thus, $p^n \mid m$.

If $\mu_{p^n}$ does not preserve $F$, then $F$ is smooth by Lemma \ref{movable} and Proposition \ref{fixedlocireducible}. The induced $\mu_p$-action does not preserve $F$ either, because otherwise it would fix $F$ pointwise by Lemma \ref{fixedloci} and then $\mu_{p^n}$ would preserve $F$ by Lemma \ref{fixedofsubscheme}. Hence, the quotient $X/\mu_{p^n}$ is smooth in a neighborhood of the image $F'$ of $F$ and the inverse image of $F'$ under the quotient map $X \to X/\mu_{p^n}$ is $p^nF$. Therefore, the multiplicity $m$ of $mF$ is divisible by $p^n$.
\qed
\vspace{5mm}

In the simpler case where $f: X \to C$ admits no multiple fibers, Theorem \ref{Main} specializes to the following corollary.

\begin{Corollary}\label{corollarymain}
Let $f: X \to C$ be an elliptic surface without multiple fibers and with $\Ker(f_*)^0$ non-trivial. Then, $\chi(X,\OO_X) = 0$ and $\Ker(f_*)^0 \in \{\mu_{p^n}, M_n, E\}$, where $n \geq 0$ and $E$ is an elliptic curve.
\end{Corollary}

\prf
Since $\Ker(f_*)^0 \cong \Ker(J(f)_*)^0$ by Corollary \ref{ExtendabilityCorollary}, we may assume that $f$ admits a section $\Sigma$. As in the proof of Theorem \ref{Main}, we have $\Sigma^2 = \chi(X,\OO_X) = 0$. In particular, $\Sigma$ satisfies $h^0(\Sigma,N_{\Sigma/X}) = 1$ and therefore $\Ker(f_*)[F] \in \{\mu_p,\alpha_p\}$ by Theorem \ref{Main}. Then, the result follows from Lemma \ref{connsub}.
\qed
\vspace{3mm}

\section{Horizontal component of $\Aut_X^0$} \label{horizontal}
Recall that the horizontal component of $\Aut_X^0$ is defined as follows.

\begin{Definition}
Let $f: X \to C$ be an elliptic surface. The \emph{horizontal component} of $\Aut_X^0$ is defined as $\im(f_*)$, where $f_*: \Aut_X^0 \to \Aut_C^0$ is the natural map.
\end{Definition}

If $f: X \to C$ is an elliptic surface such that $\im(f_*)$ is non-trivial, then certainly $\Aut_C^0$ is non-trivial and therefore $H^0(C,T_C) \neq 0$. In particular, either $C = \PP^1$ or $g(C) = 1$
Let us first treat the simpler case where the base curve $C$ satisfies $g(C) = 1$.

\begin{Proposition}\label{Propelliptic}
Let $f: X \to C$ be an elliptic surface with $g(C) = 1$. Assume that $\im(f_*)$ is non-trivial. Then, the following hold:
\begin{enumerate}[(i)]
\item All fibers of $f$ are $\Aut_X^0$-movable and $f$ is isotrivial.
\item We have $\chi(X,\OO_X) = 0$, unless possibly if $p \in \{2,3\}$, $f$ admits a multiple fiber and both the generic fiber of $f$ and $C$ are supersingular.
\item If additionally $h^0(X,T_X) \geq 2$, then one of the following holds:
\begin{enumerate}[(1)]
\item $X$ is an Abelian surface and $h^0(X,T_X) = 2$.
\item $X$ is a bielliptic surface with $\omega_X = \OO_X$ and $h^0(X,T_X) = 2$.
\item The generic fiber of $f$ is supersingular and $f$ admits a multiple fiber.
\end{enumerate}
\end{enumerate}
\end{Proposition}

\prf
First, note that if $Z \subseteq X$ is any closed subscheme contained in a fiber of $f$, then $\Aut_{(Z,X)}^0 \subseteq \Ker(f_*)$. Indeed, the action of $\Aut_{(Z,X)}^0$ on $C$ preserves the reduced point $f(Z)$ by Lemma \ref{image} and is therefore trivial by Lemma \ref{fixedloci} (ii). Hence, all fibers of $f$ are $\Aut_X^0$-movable and thus, by Lemma \ref{movable}, the fibration $f$ is isotrivial because all fibers of $f$ are either of type ${\rm II}$ or smooth and therefore the $j$-map has no poles. This proves Claim (i).

For Claim (ii), note that, by Ogg's formula, we have $\chi(X,\OO_X) = 0$ if and only if $f$ admits no fiber of type ${\rm II}$. Assume that $f$ admits a fiber $F$ of type ${\rm II}$. Then, $p \in \{2,3\}$ and the $j$-map of $f$ is identically $0$ so that the generic fiber of $f$ is supersingular and $\Aut_{(F,X)}^0 \subseteq \Ker(f_*)$ is non-trivial by Remark \ref{criteriaremark}. In particular, $f$ admits a multiple fiber by Theorem \ref{Main} (ii). Now, if $C$ is ordinary, then $\mu_p \subseteq \im(f_*)$ and since the induced extension of $\mu_p$ by $\Ker(f_*)$ splits by Lemma \ref{homsandexts} (ii) (3), there is a $\mu_p$-action on $X$. But then we can use Remark \ref{criteriaremark} again to deduce that $\mu_p \subseteq \Ker(f_*)$, which is impossible since the generic fiber of $f$ is supersingular. Thus, $C$ has to be supersingular, too.

For Claim (iii), assume that $h^0(X,T_X) \geq 2$. Let $D \in H^0(X,T_X)$ be a $p$-closed vector field. Then, we have the following short exact sequence obtained by saturating the inclusion $\OO_X \to T_X$ induced by $D$, where $Z$ (resp. $W$) is the divisorial (resp. codimension $2$) part of the zero locus of $D$:
$$
0 \to \OO_X(Z) \to T_X \to I_W(-Z-K_X) \to 0.
$$
 Note that $K_X$ is effective by the canonical bundle formula and $Z$ is effective by definition. Now, we get two cases according to whether $-Z-K_X$ is effective or not.

If $-Z-K_X$ is not effective, then $h^0(X,T_X) = h^0(X,\OO_X(Z)) \geq 2$. Therefore, the zero locus of every $p$-closed vector field contains a divisor linearly equivalent to $Z$, hence all these vector fields induce the trivial vector field on $C$. In particular, the tangent space of $\Ker(f_*)$ is at least $2$-dimensional. By Lemma \ref{connsub} and Corollary \ref{corollarymain}, this implies that the generic fiber of $f$ is supersingular and $f$ admits a multiple fiber, that is, $X$ is as in Case (3).

If $-Z-K_X$ is effective, then both $Z$ and $K_X$ are trivial. Since $h^0(X,T_X) \geq 2$, we must have $h^0(X,I_W) \geq 1$ and hence $W$ is trivial. Thus, $h^0(X,T_X) = 2$ and $X$ is Abelian or bielliptic with $\omega_X \cong \OO_X$ by the classification of surfaces. In particular, $X$ is as in Case (1) or (2).
\qed
\vspace{2mm}

Now that we understand the case where the base curve $C$ has genus $1$, it remains to treat the case where $C = \PP^1$. If $C$ is rational, then $b_1(X) \in \{0,2\}$ by \cite[Lemma 3.4]{KatsuraUeno}. In the second case we can argue in a similar fashion as in Proposition \ref{Propelliptic} because of the following lemma.

\begin{Lemma}\label{b1=2}
Let $f: X \to \PP^1$ be an elliptic surface with $b_1(X) = 2$. Then, the following hold:
\begin{enumerate}[(i)]
\item $\chi(X,\OO_X) = 0$.
\item There is an isomorphism $J(X) \cong \PP^1 \times E$, where $E$ is a general fiber of $f$.
\item The Albanese map $a_X: X \to {\rm Alb}(X)$ is a fibration over an elliptic curve and all fibers of $a_X$ are irreducible and reduced.
\item There is an action of $E$ on $X$ that induces a transitive action on ${\rm Alb}(X)$.
\end{enumerate}
\end{Lemma}

\prf
By Ogg's formula, to prove that $c_2(X) = \chi(X,\OO_X) = 0$, it suffices to show that $f$ admits no singular fibers. But this follows immediately from the criterion in\cite[Lemma 3.4]{KatsuraUeno}, since singular fibers consist of rational curves and are thus contracted by the Albanese morphism of $X$. This proves Claim (i).

Next, we prove Claim (ii) and show that $E$ acts on $X$. Since $c_2(X) = 0$, the fibration $f$ is isotrivial and $J(f)$ is a smooth elliptic fibration over $\PP^1$ by Ogg's formula. As there are no non-trivial finite \'etale covers of $\PP^1$, this implies that $J(X) \cong \PP^1 \times E$ for a general fiber $E$ of $f$. Hence, there is a finite Galois cover $C \to \PP^1$ with group $G$ such that the normalization of $X \times_{\PP^1} C$ is isomorphic to $E \times C$. Moreover, the quotient of $E \times C$ by the induced action of $G$ maps via a finite and birational map to $X$ and hence coincides with $X$. Since $X$ is smooth, the group $G$ acts via translations on the first factor of $E \times C$, for otherwise it would have an isolated fixed point and then $X$ would be singular. In particular, the translation action of $E$ on the first factor of $E \times C$ commutes with the $G$-action and thus descends to $X$.

To finish the proof, note that, by Igusa's formula, we have ${\rm rk} (\Pic(X)) \leq b_2(X) = c_2(X) + 2b_2(X) - 2 = 2$, so that all fibers of $f$ and $a_X$ have to be irreducible. Since $b_1(X) = 2$, the Albanese variety ${\rm Alb}(X)$ is an elliptic curve. Moreover, a general fiber $E$ of $f$ maps surjectively onto ${\rm Alb}(X)$, so that the target of the Stein factorization of $a_X$ is an elliptic curve, which then has to coincide with ${\rm Alb}(X)$ by the universal property of $a_X$. This also shows that the action of $E$ on ${\rm Alb}(X)$ is transitive and, because a general fiber of $a_X$ is reduced, this implies that in fact all fibers of $a_X$ are reduced.
\qed

\begin{Proposition}\label{Proprational}
Let $f: X \to \PP^1$ be an elliptic surface with $b_1(X) = 2$.  If $\im(f_*)$ is non-trivial, then the following hold:
\begin{enumerate}[(i)]
\item At most two fibers of $f$ are not $\Aut_X^0$-movable and $\chi(X,\OO_X)  = 0$.
\item If additionally $h^0(X,T_X) \geq 2$, then one of the following holds:
\begin{enumerate}[(1)]
\item $X$ is ruled over an elliptic curve.
\item $X$ is bielliptic, $\omega_X \cong \OO_X$ and $h^0(X,T_X) = 2$.
\item The generic fiber of $f$ is supersingular, $f$ admits a multiple fiber, and all fibers of $a_X$ are rational curves.
\end{enumerate}
\end{enumerate}
\end{Proposition}

\prf
The argument for Claim (i) is the same as in the proof of Proposition \ref{Propelliptic} (i), with the only difference that a subgroup scheme of $\Aut_{\PP^1}$ can have up to two fixed points by Lemma \ref{fixedloci} (i).
 Note that by Lemma \ref{b1=2} the equality $\chi(X,\OO_X) = 0$ holds even if $\im(f_*)$ is trivial. 
 
Before we start proving Claim (ii), observe that if $f$ does not admit a multiple fiber, then $X$ is ruled by the canonical bundle formula and the base curve of the ruling must be an elliptic curve, since $X$ admits an elliptic fibration. Hence, we are in Case (1) if $f$ admits a section.

Now, let us prove Claim (ii). Assume that $h^0(X,T_X) \geq 2$. Then, $\Ker((a_X)_*)$ is non-trivial. Therefore, if a general fiber $C$ of $a_X$ is smooth, then $C$ is either $\PP^1$ or an elliptic curve. Thus, in the first case, the Albanese map $a_X$ yields a ruling of $X$ over ${\rm Alb}(X)$. In the latter case, the Albanese map $a_X$ is an elliptic fibration over an elliptic curve and $a_X$ has no multiple fibers by Lemma \ref{b1=2}. Then, Corollary \ref{corollarymain} shows that the tangent space of $\Ker((a_X)_*)$ is $1$-dimensional, hence $\im((a_X)_*)$ is non-trivial. Thus, we can apply Proposition \ref{Propelliptic} (iii). It shows that $X$ is bielliptic with $\omega_X \cong \OO_X$, since an Abelian surface does not admit an elliptic fibration over $\PP^1$ and $a_X$ admits no multiple fibers by Lemma \ref{b1=2} (iii).

So, we may assume that the general fiber $C$ of $a_X$ is singular.
 Let $E$ be a general fiber of $f$. The map $E \to {\rm Alb}(X)$ factors through an \'etale morphism $A \to {\rm Alb}(X)$. Pulling back $X$ along this map, we obtain a smooth surface $X'$ with a fibration $g': X' \to A$ and an elliptic fibration $f': X' \to D$ obtained as the Stein factorization of $X' \to \PP^1$.
Now, both the action of $\Ker((a_X)_*)$ on $X$ and the action of $E$ on $X$ constructed in Lemma \ref{b1=2} lift to $X'$ and these two actions generate $\Aut_X^0$. Therefore, we have $h^0(X',T_X') \geq 2$ and $\im(f'_*)$ is non-trivial. Thus, $D$ is either an elliptic curve or $\PP^1$. If $D$ is an elliptic curve, then ${\rm Alb}(X')$ is a surface by \cite[Lemma 3.4]{KatsuraUeno} and the morphism $g': X' \to A$ factors through ${\rm Alb}(X')$. This is impossible, since the fibers of $g'$ are singular whereas the fibers of ${\rm Alb}(X')$ are (unions of) elliptic curves. Hence, we must have $D \cong \PP^1$ and we may replace $X$ by $X'$ to assume that a general fiber of $f$ maps purely inseparably to ${\rm Alb}(X')$ and in particular all multiple fibers of $f$ have multiplicity $p^n$ for some $n \geq 1$. The remainder of the proof splits into two cases according to whether the generic fiber of $f$ is ordinary or supersingular.

If the generic fiber of $f$ is ordinary then so is ${\rm Alb}(X)$. Let $m_iF_i$ be a multiple fiber of $f$. Since the map $F_i \to X \to {\rm Alb}(X)$ is purely inseparable, the dual map $\Pic^0_{{\rm Alb}(X)} \to \Pic^0_X \to \Pic^0_{F_i}$ is \'etale. In particular, the map $H^1({\rm Alb}(X), \OO_{{\rm Alb}(X)}) \to H^1(X,\OO_X) \to H^1(F_i,\OO_{F_i})$ is an isomorphism. By \cite[Section 6]{KatsuraUeno}, one can use cocycles $\rho \in H^1(X,\OO_X)$ that are fixed by Frobenius and map to a non-trivial element in $H^1(F_i,\OO_{F_i})$ to construct an \'etale cover $\widetilde{X}$ of $X$ with an elliptic fibration $\widetilde{f}: \widetilde{X} \to D$ without multiple fibers. Choosing the cocycles $\rho$ in the image of $H^1({\rm Alb}(X), \OO_{{\rm Alb}(X)}) \to H^1(X,\OO_X)$, we can assume that $\widetilde{X}$ arises as pullback of $a_X$ along an \'etale isogeny $A \to {\rm Alb}(X)$. Then, as in the previous paragraph, the group scheme $\Aut_X^0$ acts on $\widetilde{X}$ and $\im(\widetilde{f}_*)$ is non-trivial and the image $D$ of the Stein factorization of $\widetilde{X} \to \PP^1$ satisfies $D \cong \PP^1$. Since $\widetilde{f}$ admits no multiple fibers and $\PP^1$ admits no \'etale covers, the equality $c_2(\widetilde{X}) = c_2(X) = 0$ implies $\widetilde{X} \cong \PP^1 \times A$, contradicting our assumption that the fibers of $a_X$ are singular.

Hence, the generic fiber of $f$ is supersingular. Then, it is shown in \cite[Proposition 3.1]{Kawazoe} that there is a purely inseparable cover $\widetilde{\pi}: \widetilde{X} \to X$ such that the Stein factorization of $f \circ \widetilde{\pi}$ is an elliptic fibration $\widetilde{f}: \widetilde{X} \to D$ without multiple fibers. Since $\widetilde{\pi}$ is purely inseparable, so is $D \to \PP^1$. Hence, $D \cong \PP^1$ and $\widetilde{f}$ admits a section $\Sigma$. The image $\widetilde{\pi}(\Sigma)$ is a rational curve and is therefore contracted by $a_X$. The fibers of $a_X$ are integral, so that $\widetilde{\pi}(\Sigma)$ coincides with a fiber of $a_X$. Hence, all fibers of $a_X$ are rational curves.
\qed
\vspace{5mm}

Thus, the last remaining case are elliptic surfaces $f: X \to \PP^1$ with $b_1(X) = 0$. We will use the following lemma, which is well-known in characteristic $0$.

\begin{Lemma}\label{twofibers}
Let $f: X \to \PP^1$ be an elliptic surface with $b_1(X) = 0$ and at most two singular fibers $F_1,F_2$ with Swan conductors $\delta_{F_1} = \delta_{F_2} = 0$. Then, $\chi(X,\OO_X) = 1$ and the possible types of $F_1$ and $F_2$ are as follows:
\begin{enumerate}[(i)]
\item $({\rm II},{\rm II}^*)$ and $p \not\in \{ 2,3 \}$.
\item $({\rm III},{\rm III}^*)$ and $p \neq 2$.
\item $({\rm IV},{\rm IV}^*)$ and $p \neq 3$.
\item $({\rm I}_0^*,{\rm I}_0^*)$ and $p \neq 2$.
\end{enumerate}
\end{Lemma}

\prf
Since $f$ and $J(f)$ have the same types of singular fibers and $\chi(X,\OO_X) = \chi(J(X),\OO_{J(X)})$, we may assume that $f$ admits a section. Let $F_1,F_2$ be the singular fibers of $f$.

Let $T = U \oplus T_1 \oplus T_2$, where $U$ is the unimodular lattice generated by a section of $f$ and the class of a fiber of $f$, and $T_i$ is the lattice generated by the components of $F_i$ disjoint from the zero section of $J(f)$. Then, Igusa's inequality yields
$$
{\rm rk}(T) = 2 + {\rm rk}(T_1) + {\rm rk}(T_2)  \leq {\rm rk}(\Num(X)) \leq b_2(X) = c_2(X) - 2.
$$
On the other hand, by Ogg's formula and our assumption that the Swan conductor of every fiber is trivial, we have
$$
c_2(X) - 2 = e(F_1) + e(F_2) - 2 \leq 2 + {\rm rk}(T_1) + {\rm rk}(T_2)
$$
with equality if and only if $F_1$ and $F_2$ are additive. Thus, both $F_1$ and $F_2$ are additive fibers and $T \subseteq \Num(X)$ is of finite index. Moreover, $T_1$ and $T_2$ are either trivial or root lattices of type $A_1,A_2,D_n,E_6,E_7$ or $E_8$. Their discriminants are $2,3,4,3,2$ and $1$, respectively.

Since ${\rm rk} (\Num(X)) = b_2(X)$, we can use $\ell$-adic Poincar\'e duality for all $\ell \neq p$ to deduce that the discriminant ${\rm disc}(\Num(X))$ is a power of $p$. Moreover, since $T \subseteq \Num(X)$ is of finite index, the discriminants of these two lattices differ by a square. Taking into account ${\rm rk}(T) = b_2(X) = c_2(X) - 2 = 10 +12k$ for some $k \geq 0$, we thus have the following cases, where in each case we have ${\rm disc}(\Num(X)) = 1$:
\begin{enumerate}[(i)]
\item $T_1 = 0$, $T_2 = E_8$ and $p \neq 2,3$,
\item $T_1 = A_1$, $T_2 = E_7$ and $p \neq 2$,
\item $T_1 = A_2$, $T_2 = E_6$ and $p \neq 3$,
\item $T_1 = D_m$, $T_2 = D_n$ for some $m,n$ and $p \neq 2$,
\item $T_1 = D_m$, $T_2 \in \{E_8,0\}$ for some $m$ and $p \neq 2$.
\end{enumerate}
Now, if $p \neq 2$, we can apply a quadratic twist to $f$ that only changes the fibers $F_1$ and $F_2$. Then, either all fibers of the twisted fibration are smooth or the fibration satisfies the assumptions of the lemma and then its singular fibers have to appear in the above list. Hence, Lemma \ref{twist} (i) shows that $m = n = 4$ in Case $(iv)$ and that Case $(v)$ does not exist.
\qed
\vspace{2mm}

\begin{Remark} \label{twofibersremark}
The Jacobian $J(X)$ of each of these four types of surfaces in the above Lemma \ref{twofibers} is a rational surface that can be defined over $\ZZ[j(F_\eta)]$, where $F_\eta$ is the generic fiber of $f$. The reductions of these fibrations modulo the excluded characteristics in the respective cases in Lemma \ref{twofibers} become quasi-elliptic (see e.g. \cite{Lang3}). It is straightforward, e.g. from the equations given in \cite{Lang3}, to check that $J(X)$ is the minimal resolution of $(\PP^1 \times E)/(\ZZ/n\ZZ)$ with $n \in \{2,3,4,6\}$ and $\ZZ/n\ZZ$ acting diagonally with a fixed point on $E$.
\end{Remark}

If we allow the Swan conductors $\delta_{F_i}$ to be non-trivial in Lemma \ref{twofibers}, there are many examples of elliptic surfaces with only one or two singular fibers. In the following lemma, we will treat a very special case that will appear in Theorem \ref{mainrational}.

\begin{Lemma} \label{twofibers23}
Let $f: X \to \PP^1$ be an elliptic surface with at most two singular fibers $F_1,F_2$. Assume that $\delta_{F_2} = 0$, and either $p = 3$ and $F_1$ is of type ${\rm II}$ with $\delta_{F_1} = 1$ or $p = 2$ and $F_1$ is of type ${\rm II}$ or ${\rm III}$ with $\delta_{F_1} = 2$ or $\delta_{F_1} = 1$, respectively. Then, the following hold:
\begin{enumerate}[(i)]
\item If $p = 3$, the possible types of $F_1$ and $F_2$ are $({\rm II},{\rm III}^*),({\rm II},{\rm I}_{3^{2k}}),({\rm II},{\rm I}_{3^{2k-1}}^*)$. In particular, $c_2(X) = 3^{2k} + 3$ or $c_2(X) = 3^{2k-1} + 9$ for some $k \geq 1$.
\item If $p = 2$, the possible types of $F_1$ and $F_2$ are $({\rm II},{\rm IV}^*),({\rm III},{\rm IV}^*),({\rm II},{\rm I}_{2^{2k+1}}),({\rm III},{\rm I}_{2^{2k+1}})$. In particular, $c_2(X) = 2^{2k+1} + 4$ for some $k \geq 1$.
\end{enumerate}
\end{Lemma}

\prf
 As in the proof of Lemma \ref{twofibers}, we may assume that $f$ admits a section. We will split the proof in two cases according to whether $f$ is isotrivial or not.

Assume first that $f$ is isotrivial. Then, since the $j$-map of $f$ has a zero at $f(F_1)$, the generic fiber of $f$ is supersingular with $j$-invariant $0$. If $p = 3$, then $v_{f(F_1)}(\Delta_f) = 3$ by assumption so that $v_{f(F_2)}(\Delta_f) \equiv 9$ mod $12$ by Ogg's formula. Moreover, we have $\delta_{F_2} = 0$ and the $j$-map has no pole at $f(F_2)$, so that $F_2$ is additive with $8$ components. This implies that $F_2$ is of type ${\rm III}^*$. If $p = 2$, then $v_{f(F_1)}(\Delta_f) = 4$ and thus, by the same argument as before, $F_2$ is additive with $7$ components. Since $\delta_{F_2} = 0$, this implies that $F_2$ is of type ${\rm IV}^*$.

Next, assume that $f$ is not isotrivial. Then, the $j$-map of $f$ has a pole. If $p = 3$, this implies that $F_2$ is of type ${\rm I}_n^*$ or ${\rm I}_n$. Then, by Lemma \ref{twist}, we can replace $f$ by a quadratic twist and assume that $F_2$ is of type ${\rm I}_n$ and $F_1$ is of type ${\rm II}$ or ${\rm II}^*$ with $\delta_{F_1} = 1$.
If $p = 2$, the assumption $\delta_{F_2} = 0$ forces $F_2$ to be of type ${\rm I}_n$. Moreover, by Lemma \ref{twist}, we can replace $f$ by a quadratic twist to assume that $F_1$ is of type ${\rm III}$ with $\delta_{F_1} = 1$.

Now, we let $T = U \oplus T_1 \oplus T_2$, where $T_i$ is spanned by non-identity components of $F_i$, and $U$ is generated by the class of a fiber and a section of $f$. Note that $T_1$ is unimodular by the previous paragraph. Since $F_2$ is multiplicative and $F_1$ satisfies $\delta_{F_1} = 1$, we obtain ${\rm rk} (T) = b_2(X)$ from Ogg's formula. Hence, $T$ has finite index in $\Num(X)$. As in the proof of Lemma \ref{twofibers}, $\ell$-adic Poincar\'e duality shows that ${\rm disc}(\Num(X))$ is a power of $p$. By \cite[Section 11.10]{schuettShioda} this implies that $n = p^i m^2$, where $i$ is some integer and $m$ is the order of the group of torsion sections of $f$ of order prime to $p$. Since, on the one hand, a torsion section of order prime to $p$ is disjoint from the zero section \cite[Proposition 3.5 (iv)]{OguisoShioda} and, on the other hand, $f$ admits the fiber $F_1$ whose underlying group is $\GG_a$, we have $m = 1$. Thus, we have $v_{f(F_1)}(\Delta_f) + p^i \equiv 0$ mod $12$ by Ogg's formula.

If $p = 3$ and $F_1$ is of type ${\rm II}$, this implies that $i = 2k$ for some $k \geq 1$, and if $F_1$ is of type ${\rm II}^*$, this implies that $i = 2k-1$ for some $k \geq 1$. Undoing the quadratic twist we applied in the second paragraph of the proof, we obtain the stated types of singular fibers.

If $p = 2$, then this implies $i = 2k+1$ for some $k \geq 1$. Again, undoing the quadratic twist, we obtain the stated types of singular fibers.
\qed
\vspace{2mm}

After having prepared the necessary technical lemmas, we are now ready to prove the main result of this section on elliptic surfaces $f: X \to \PP^1$ with $b_1(X) = 0$.

\begin{Theorem} \label{mainrational}
Let $f: X \to \PP^1$ be an elliptic surface with $b_1(X) = 0$. Assume that $\im(f_*)$ is non-trivial. Then, $f$ has at most two non-$\Aut_X^0$-movable fibers and one of following holds:
\begin{enumerate}[(i)]
\item $f$ is isotrivial with precisely two singular fibers of the types given in Lemma \ref{twofibers}. Moreover, $\Aut_X^0 \subseteq \GG_m$.
\item $p \in \{2,3\}$ and $f$ admits precisely two singular fibers $F_1,F_2$ of the types given in Lemma \ref{twofibers23}. Moreover, $\Aut_X^0 \cong \im(f_*) \cong \mu_p$ and there are no multiple fibers except possibly $F_1$ and $F_2$.
\item $p = 2$, the generic fiber of $f$ is ordinary, $f$ admits a fiber $F$ of type ${\rm I}_{8k + 4}^*$ with $\delta_{F} = 4k+8$ for some $k \geq 0$ and all other fibers of $f$ are smooth and non-multiple. Moreover, we have $\Aut_X^0 \cong \im(f_*) \subseteq \GG_a$.
\item $p \in \{2,3\}$ and $f$ is isotrivial with supersingular generic fiber and at most one non-$\Aut_X^0$-movable fiber. Moreover, the group scheme $\Aut_X^0$ does not contain $\mu_p$.
\end{enumerate}
\end{Theorem}

\prf
Since $b_1(X) =  0$, we have $c_2(X) > 0$ and thus $f$ admits at least one singular fiber $m_1F_1$, say over $\infty \in \PP^1$, where $m_1$ is the multiplicity of the fiber. On the other hand, by the same argument as in Proposition \ref{Proprational}, $f$ admits at most two non-$\Aut_X^0$-movable fibers. To prove the remaining claims, we will make use of fact that we understand the fixed loci of $\mu_p$-actions on $X$ by Theorem \ref{fixedpointformula}. To do this, we will split the proof into three cases according to whether $\mu_p \subseteq \Ker(f_*)$, $\mu_p \subseteq \im(f_*)$ or $\Aut_X^0$ does not contain any $\mu_p$ at all.

\vspace{2mm}
\underline{Case $\mu_p \subseteq \Ker(f_*)$:}

Assume that $\mu_p \subseteq \Ker(f_*)$. Then, by Lemma \ref{connsub}, we have $\mu_p = \Ker(f_*)[F]$ and in particular $\mu_p$ is preserved by every automorphism of $\Ker(f_*)$. Since $\Ker(f_*)$ is normal in $\Aut_X^0$, this implies that $\mu_p$ is normal in $\Aut_X^0$. Therefore, Lemma \ref{fixedofsubscheme} implies that the action of $\Aut_X^0$ on $X$ preserves $X^{\mu_p}$. Since $\im(f_*)$ is non-trivial by assumption, we can apply Lemma \ref{fixedloci} to deduce that there is a fiber $m_2F_2$ of $f$, say over $0 \in \PP^1$, such that $X^{\mu_p} \subseteq F_1 \cup F_2$, for otherwise the action of $\Aut_X^0$ on $\PP^1$ would have more than two fixed points, which is impossible. By  Remark \ref{criteriaremark} and Lemma \ref{fixedloci}, the $\mu_p$-action on $X$ preserves every singular fiber of $f$ and has at least one fixed point on each such fiber. In particular, $m_1F_1$ and $m_2F_2$ are the only possibly singular fibers of $f$. The fixed point formula given in Theorem \ref{fixedpointformula} then yields
\begin{equation}\label{fixedpointeqn}
c_2(X) = e(F_1^{\mu_p}) + e(F_2^{\mu_p}) \tag{$\ast$}
\end{equation}

Assume that $p \not \in \{2,3\}$. Then, Lemma \ref{twofibers} implies that both $F_1$ and $F_2$ are singular of the types given in Lemma \ref{twofibers}.

If $p \in \{2,3\}$, we can compare Equation \eqref{fixedpointeqn} with Ogg's formula to obtain
$$
e(F_1^{\mu_p}) + e(F_2^{\mu_p}) = c_2(X) = e(F_1) + e(F_2) + \delta_{F_1} + \delta_{F_2}.
$$
By Proposition \ref{fixedlocireducible}, we know that $e(F_i^{\mu_{p}}) = e(F_i)$ holds unless $F_i$ is of type ${\rm II}$, or $F_i$ is of type ${\rm III}$ and $p = 2$. If neither $F_1$ nor $F_2$ are of these types, then $\delta_{F_1} = \delta_{F_2} = 0$ and we conclude as in the case $p \not \in \{2,3\}$. Note also that if $F_i$ is of type ${\rm II}$ or ${\rm III}$, then $e(F_i^{\mu_p}) \leq 4$ by Lemma \ref{fixedloci} and Proposition \ref{fixedlocireducible}, so that $c_2(X) = 12 \chi(X,\OO_X) \geq 12$ implies that not both $F_1$ and $F_2$ are of these exceptional types. Hence, we may assume that $e(F_2) = e(F_2^{\mu_p})$.

Assume that $p = 3$ and $F_1$ is of type ${\rm II}$. Then, we have $e(F_1^{\mu_3}) = 3 = e(F_1) + 1$ by Proposition \ref{fixedlocireducible}, and hence Equation \eqref{fixedpointeqn} shows that $\delta_{F_1} + \delta_{F_2} = 1$. Since $\delta_{F_1} \geq 1$, this implies $\delta_{F_2} = 0$ and thus we can apply Lemma \ref{twofibers23} to determine the types of $F_1$ and $F_2$. 

If $p = 2$ and $F_1$ is of type ${\rm III}$, then $e(F_1^{\mu_2}) = 4 = e(F_1) + 1$ and Equation \eqref{fixedpointeqn} shows that $\delta_{F_1} + \delta_{F_2} = 1$. The rest of the argument is as in the case $p = 3$. Similarly, if $F_1$ is of type ${\rm II}$, then $\delta_{F_1} \geq 2$, so that again $\delta_{F_2} = 0$ and Lemma \ref{twofibers23} applies.

We have shown that the singular fibers of $f$ are as claimed in (i), (ii), or (iii) and it remains to prove the assertions on the structure of $\Aut_X^0$ and the multiple fibers. For this, we will first show that $h^0(X,T_X) \leq 1$ holds.
Denote the divisorial part of $X^{\mu_p}$ by $Z$ and the isolated part by $W$. Then, the saturation of the section of $T_X$ given by the $\mu_p$-action yields an exact sequence
$$
0 \to \OO_X(Z) \to T_X \to I_W(- K_X - Z) \to 0.
$$
Since $F_1$ and $F_2$ are singular and $Z$ is smooth, the $F_i$ cannot be contained in $Z$. Hence, we have $h^0(X,\OO_X(Z)) \leq 1$ and the above sequence shows that $h^0(X,T_X) \leq 1$, unless possibly if $- K_X$ is effective. If $-K_X$ is effective, then the canonical bundle formula shows that $f$ admits no wild fibers and at most one multiple fiber. In this case, if $f$ admits no multiple fiber, then $\Ker(f_*)$ is trivial by Theorem \ref{Main} so that $h^0(X,T_X) \leq 1$. If $f$ admits a multiple fiber $F$, then $h^0(X,\OO_X(-K_X)) = 1$ and $-K_X \sim F$. In this case, we also have $h^0(X,T_X) \leq 1$, unless $X^{\mu_p} \subseteq F$. But Lemma \ref{fixedloci} shows that $\mu_p$ has fixed points on both $F_1$ and $F_2$, so $X^{\mu_p} \subseteq F$ is impossible. Therefore, we have $h^0(X,T_X) \leq 1$ in all cases.

Now, since $\mu_p$ has fixed points on $F_1$ and $F_2$ and $\Aut_X^0$ acts on $X^{\mu_p}$, we have $\im(f_*) \subseteq \Aut^0_{0 \cup \infty, \PP^1} \cong \GG_m$. By Theorem \ref{Main}, we also have $\Ker(f_*) \cong \mu_{p^n}$ for some $n \geq 0$. Thus, the group scheme $\Aut_X^0[F^n]$, being an extension of finite commutative group schemes of type $lr$, is also commutative of type $lr$ by Lemma \ref{homsandexts}. Since $h^0(X,T_X) \leq 1$, Lemma \ref{finitecommgroupscheme} implies $\Aut_X^0[F^n] \cong \mu_{p^n}$ and therefore either $\Aut_X^0 \cong \GG_m$ or $\Aut_X^0 \cong \mu_{p^n}$ for some $n \geq 1$. In the cases where $p \in \{2,3\}$ and $F_1$ is of type ${\rm II}$ or ${\rm III}$, Lemma \ref{fixedloci} and Proposition \ref{fixedlocireducible} imply that $n = 1$. But then $\Aut_X^0$ acts trivially on the base, contradicting $\mu_p \subseteq \Ker(f_*)$. Putting everything together, we see that $\mu_p \subseteq \Ker(f_*)$ implies that we are in Case (i).

\vspace{2mm}
\underline{Case $\mu_p \not \subseteq \Ker(f_*)$ and $\mu_p \subseteq \im(f_*)$:}

If $\mu_p \subseteq \im(f_*)$ and $\mu_p \not \subseteq \Ker(f_*)$, then by Lemma \ref{connsub} $\Ker(f_*)$ is either trivial or finite and commutative of type $ll$. Thus, the extension of $\mu_p$ by $\Ker(f_*)$ splits by Lemma \ref{homsandexts} and we get a $\mu_p$-action on $X$ whose fixed locus is contained in two fibers, which are then necessarily the only singular or multiple fibers of $f$. Then, the arguments where we compute the types of $F_1$ and $F_2$ and deduce $h^0(X,T_X) \leq 1$ are the same as in the previous case. But this shows that $\Ker(f_*)$ has to be trivial, for otherwise $\Aut_X^0$ would contain $\alpha_p \rtimes \mu_p$ and thus its tangent space would be too big. Moreover, the fiber $F_2$ is not $\Aut_X^0$-movable by Lemma \ref{movable}, so $\Aut_X^0 \cong \im(f_*) \subseteq \Aut^0_{(0,\PP^1)} \cong \GG_a \rtimes \GG_m$. The only subgroup schemes of $\GG_a \rtimes \GG_m$ which have $1$-dimensional tangent space and contain $\mu_p$ are $\mu_{p^n}$ and $\GG_m$, so $\Aut_X^0$ has to be one of those two group schemes. Moreover, in the cases where $p \in \{2,3\}$ and $F_1$ is of type ${\rm II}$ or ${\rm III}$, Lemma \ref{fixedloci} and Proposition \ref{fixedlocireducible} imply that $\Aut_X^0 \cong \mu_p$. Thus, $\mu_p \not \subseteq \Ker(f_*)$ and $\mu_p \subseteq \im(f_*)$ imply that we are in Case (i) or (ii).

\vspace{2mm}
\underline{Case $\mu_p \not \subseteq \Ker(f_*)$ and $\mu_p \not \subseteq \im(f_*)$:}

Since $\mu_p \not \subseteq \im(f_*)$, Lemma \ref{fixedloci} shows that the action of $\Aut_X^0$ on $\PP^1$ has at most one fixed point and thus $f$ has at most one non-$\Aut_X^0$-movable fiber. Hence, by Lemma \ref{movable} and Lemma \ref{twofibers}, we have $p \in \{2,3\}$.
If the generic fiber of $f$ is supersingular, then we are in Case (iv), so we may assume that the generic fiber of $f$ is ordinary.

Assume that $f$ has ordinary generic fiber. Since $\mu_p \not \subseteq \Ker(f_*)$, we have $\Ker(f_*) = \{ {\rm id}\}$. The $j$-map is not identically $0$, so $F_1$ is not of type ${\rm II}$ and in particular not $\Aut_X^0$-movable by Lemma \ref{movable}. Hence, $\Aut_X^0$ acts on $\PP^1$ with a fixed point and thus $\Aut_X^0 \subseteq \Aut_{(\infty,\PP^1)}^0 \cong \GG_a \rtimes \GG_m$. In fact, by our assumption that $\mu_p \not \subseteq \im(f_*)$, we have $\Aut_X^0 \subseteq \GG_a$.
In particular, there is an $\alpha_p \subseteq \Aut_X^0$ that acts non-trivially on $\PP^1$. By Lemma \ref{movable} singular fibers are preserved by $\alpha_p$, hence $F_1$ is the only singular fiber of $f$. Similarly, if $f$ admits a multiple fiber $mF$ different from $F_1$, then $mF$ is $\Aut_X^0$-movable so that $pF \subseteq mF$ by Lemma \ref{movable}. But then $\alpha_p$ preserves $pF$ by Lemma \ref{Frobeniuspower}, contradicting the fact that $\alpha_p$ acts with only one fixed point on $\PP^1$. Therefore, all multiple or singular fibers of $f$ are equal to $F_1$.
By Lemma \ref{twist} (iii), this implies that $f$ is isotrivial and $F_1$ is of type ${\rm I}_{8k+4}^*$ with $\delta_{F_1} = 4k + 8$ for some $k \geq 0$. Hence, we are in Case (iii). This finishes the proof.
\qed
\vspace{15mm}

\section{Examples}\label{examples}
The purpose of this section is to give several examples illustrating the different phenomena discussed in the previous two sections.

\subsection{Examples with many global vector fields} \label{Vectorfieldssection}
In this section, we show that all types of surfaces $X$ with $h^0(X,T_X) \geq 2$ listed in Theorem (A) actually occur. Moreover, we give a series of examples proving Theorem (A) (iii). Throughout, $E$ denotes an elliptic curve.

\begin{Example}[\emph{Elliptic ruled surfaces}]\label{ruled}
If $X$ is ruled over $E$, let $a_X: X \to E$ be the ruling and assume that $X$ admits an elliptic fibration $f: X \to \PP^1$. Being a ruled surface, $X$ can be written as $X = \PP(\cal{E})$ for some normalized (in the sense of \cite[Chapter V, Proposition 2.8]{Hartshorne}) locally free sheaf $\cal{E}$ of rank $2$ on $E$. Let $e := -\deg(\cal{E})$. Using the results of \cite[Chapter V, Corollary 2.18]{Hartshorne}, it is straightforward to check that $e \in \{0,-1\}$. Therefore, either $\cal{E}$ is the unique indecomposable vector bundle of rank $2$ on $E$ with $e \in \{0,-1\}$ or $\cal{E} \cong \OO_E \oplus \cal{L}$ for a torsion line bundle $\cal{L}$ on $E$ of order $n \geq 0$. Finally, it follows from \cite[Theorem 9]{Atiyah} that if $p = 0$ and $\cal{E}$ is indecomposable with $e = 0$, then $X$ does not admit an elliptic fibration while \cite[Proposition, p.336]{Mumford} implies that the corresponding ruled surface admits an elliptic fibration in positive characteristic.

Alternatively, these surfaces can be described as $X = (E \times \PP^1)/G$, where $G \subseteq E$ is a finite subgroup scheme acting faithfully on $\PP^1$. With this description, it is clear that if $N$ is the normalizer of $G$ in $\Aut_{E \times \PP^1}$, then $\Aut_X = N/G$. Since $\Aut_{E \times \PP^1}^0 \cong E \times {\rm PGL}_2$ and $E$ is commutative, we can calculate $N^0$ as the product of the centralizers of $G$ in $E$ and ${\rm PGL}_2$, respectively. Putting all of this together, the connected component of the automorphism scheme of an elliptic surface which is also ruled over an elliptic curve $E$ is as in the following table:

\begin{table}[h!]
\resizebox{\textwidth}{!}{$\displaystyle
\begin{array}{|l|l|l|l|} \hline
\cal{E} & G & \Aut_X^0/E & h^0(X,T_X) \\ \hline
\OO_E \oplus \OO_E & \{1\} & {\rm PGL}_2 & 4 \\ \hline
\OO_E \oplus \cal{L} & \mu_n & \GG_m & 2 \\ \hline
indec., e = -1 & 
E[2] & \{1\} & 1 \\ \hline
indec., e = 0 & \begin{cases}
\ZZ/p\ZZ &\text{ if } E \text{ is ordinary} \\
\alpha_p  &\text{ if } E \text{ is supersingular}
\end{cases} & 
\begin{cases}
\GG_a &\text{ if } p \neq 2\\
\GG_a \rtimes \mu_2 & \text{ if } p = 2 \text{ and } E \text{ is ordinary} \\
\GG_a \times \alpha_2 & \text{ if } p = 2 \text{ and } E \text{ is supersingular}
\end{cases} & \begin{cases}
2 &\text{ if } p \neq 2\\
3 &\text{ if } p = 2
\end{cases} \\ \hline
\end{array} $
}
\end{table}
\noindent 
The calculation of $(\Aut_X^0)_{red}$ and $h^0(X,T_X)$ for all ruled surfaces can be found in \cite{Maruyama}. Therefore, the only thing in the above table that still needs to be checked is the case $p = 2$ and $e = 0$ and we leave this case to the reader.
\end{Example}

\begin{Example}[\emph{Abelian and bielliptic surfaces}]\label{abelianbielliptic}
If $X$ is Abelian, then $\Aut_X^0 \cong X$ and in particular $h^0(X,T_X) = 2$. If $X$ is bielliptic, then, by \cite{BombieriMumford3}, the canonical sheaf $\omega_X$ can be trivial if and only if $p \in \{2,3\}$. In these cases, one can prove that $\Aut_X^0$ is not reduced. We refer the reader to the article \cite{Martin} of the author, where the group scheme $\Aut_X$ is calculated for all (quasi-)bielliptic surfaces in all characteristics.
\end{Example}

\begin{Example}[\emph{Examples with supersingular generic fiber}]\label{supersingularexample}
This example will serve as the basic example of isotrivial elliptic surfaces with supersingular generic fiber and many vector fields from which we will derive a whole series of examples in Example \ref{importantexample}.
Consider the rational curve $C \subseteq \PP^2$ of degree $p+1$ given by the homogeneous equation
$$
y^pz = x^{p+1}.
$$
Then, $p_a(C) = \frac{p(p-1)}{2}$ and $C$ has a single isolated singularity at $P = [0:0:1]$. Consider the $\alpha_p$-action defined by
$$
[x:y:z] \mapsto [x: y+az:z] \hspace{1cm} a^p = 0
$$
and note that $P$ is not a fixed point of this action. More precisely, the reduced fixed locus of $\alpha_p$ on $C$ consists of the single smooth point $Q = [0:1:0]$.

Now, let $E$ be a supersingular elliptic curve and let $X := (E \times C)/\alpha_p$, where $\alpha_p \subseteq E$ acts on $C$ via the action defined above. By the same argument as in the proof of \cite[Proposition 7]{BombieriMumford3}, the surface $X$ is smooth, since $\alpha_p$ does not fix $P$. Moreover, $X$ comes with two fibrations $a_X: X \to E/\alpha_p$ and $f: X \to \PP^1$, where the latter is obtained by taking the normalization of $C/\alpha_p$. By construction, the morphism $f$ is an elliptic fibration with general fiber isomorphic to $E$ and $f$ admits a unique multiple fiber of multiplicity $p$, namely the image of $E \times Q$ on $X$.

Finally, note that there is an $\alpha_p^2$-action on $C$ given by
$$
[x:y:z] \mapsto [x^2 + bxy: xy + cx^2 + bcxy: xz + byz] \hspace{1cm} b^p = c^p = 0
$$
and this action commutes with the $\alpha_p$-action used to construct $X$. Thus, we get an induced action of $E \times \alpha_p^2$ on $X$. In particular, we have $h^0(X,T_X) \geq 3$.

\begin{Remark}
Alternatively, and analogously to the construction of Raynaud's counterexamples to Kodaira vanishing on surfaces in positive characteristic \cite{Raynaud}, one can describe the above example as follows: Let $E$ be a supersingular elliptic curve and let $\cal{E}$ be the indecomposable vector bundle of rank $2$ on $E$ with $e = 0$. Since $E$ is supersingular, the Frobenius map $F$ is trivial on $H^1(E,\OO_E) \cong {\rm Ext}^1(\OO_X,\OO_X)$, so the pullback of $\cal{E}$ along $F$ splits and this splitting yields an inseparable multisection of the ruling $\PP(\cal{E}) \to E$. Then, $X$ can be defined as the degree $(p+1)$ cover of $\PP(\cal{E})$ branched over the inseparable multisection and a disjoint section.
\end{Remark}
\end{Example}

\begin{Example}[\emph{Examples with unbounded vector fields}]\label{importantexample}
Here, we will use Example \ref{supersingularexample} to construct the elliptic surfaces announced in Theorem (A) (iii). More precisely, for every field $K$ of positive characteristic, for every smooth projective curve $\widetilde{C}$ over $K$ and for every $n \geq 1$, we will use $f: X \to C$ to construct an elliptic surface $\widetilde{f}:\widetilde{X} \to \widetilde{C}$ with $h^0(\widetilde{X},T_{\widetilde{X}}) \geq n$.

Let $E$ be a supersingular elliptic curve over $\FF_p$. Then, the surface $X := (E \times C)/\alpha_p$ constructed in Example \ref{supersingularexample} is also defined over $\FF_p$. Moreover, the elliptic fibration $f: X \to \PP^1$ has exactly one multiple fiber, corresponding to the unique fixed point of the $\alpha_p$-action on $C$.

Now, let $n \geq 1$ be arbitrary, let $K$ be some field extension of $\FF_p$ and let $\widetilde{C}$ be a smooth projective curve over $K$. If $K$ is finite, choose a finite separable morphism $g': \widetilde{C} \to \PP^1$ which is ramified over only one point (this is possible by the "wild Belyi Theorem" \cite[Theorem 1]{AnbarTutdere}) and let $g$ be the composition of $g'$ with a tame finite map $\PP^1 \to \PP^1$ of degree at least $n$. If $K$ is infinite, let $g$ be any finite and separable map $g: \widetilde{C} \to \PP^1$. In both cases, we can modify $g$ by an automorphism of $\PP^1$ such that the multiple fiber of $f$ does not map to a branch point of $g$.

Now, let $\widetilde{f}: \widetilde{X} \to \widetilde{C}$ be the base change of $f$ along $g$. The branch locus of $\widetilde{X} \to X$ consists of a disjoint union of smooth fibers, so $\widetilde{X}$ is smooth. We claim that $h^0({\widetilde{X}},T_{\widetilde{X}}) \geq n$. By flat base change, we may assume that $K$ is algebraically closed. Since $\Ker(f_*)$ preserves the fibers of $f$, it acts naturally on the fiber product $\widetilde{X} = \widetilde{X} \times_{\PP^1} \widetilde{C}$ and we obtain an inclusion $\Ker(f_*) \subseteq \Ker(\widetilde{f}_*)$. Next, consider the short exact sequence associated to any $D \in H^0(X,T_X)$, where $Z$ is the divisorial part and $W$ is the isolated part of the zero locus of $D$
$$
0 \to \OO_X(Z) \to T_X \to I_W(-Z-K_X) \to 0.
$$
Since $H^0(X,T_X)$ contains a $3$-dimensional subspace generated by the additive vector fields corresponding to the $\alpha_p^3$-action constructed in the previous example, we must have $h^0(X,\OO_X(Z)) \geq 2$ and there is an $\alpha_p$-action $\rho$ on $X$ that fixes a simple fiber $F$ of $f$ such that $f(F)$ is not a branch point of $g$. By construction, the preimage of $F$ in $\widetilde{X}$ consists of at least $n$ disjoint simple fibers $\widetilde{F}_1,\hdots,\widetilde{F}_n$, all of which must be fixed pointwise by the $\alpha_p$-action $\widetilde{\rho}$ on $\widetilde{X}$ that induces the action $\rho$. Then, we consider the short exact sequence induced by the action $\widetilde{\rho}$, where $\widetilde{Z}$ and $\widetilde{W}$ are the divisorial and isolated part of the fixed locus, respectively:
$$
0 \to \OO_{\widetilde{X}}(\widetilde{Z}) \to T_{\widetilde{X}} \to I_{\widetilde{W}}(-\widetilde{Z}-K_{\widetilde{X}}) \to 0.
$$
Since $\bigcup_{i=1}^n F_i \subseteq \widetilde{Z}$, we have $n \leq h^0(X,\OO_X(\sum_{i=1}^n F_i)) \leq h^0(X,T_X)$ by a Clifford argument. In particular, for every curve $C$ over every field $K$ of positive characteristic, the set of numbers
$$
\{ h^0(X,T_X) \mid X \text{ admits an elliptic fibration } f:X \to C \}
$$
is unbounded.
\end{Example}

\subsection{Non-isotrivial examples with non-trivial vertical component and $\Aut_X^0$-movable multiple fibers}
In this section, we let $k$ be an algebraically closed field of characteristic $2$. We give examples of non-isotrivial elliptic surfaces $f: X \to C$ such that $\Ker(f_*)$ is non-trivial and also show that $\Aut_X^0$-$1$-movable multiple fibers exist over $k$, even for non-isotrivial elliptic surfaces.
Recall that a supersingular Enriques surface $X$ over $k$ is a smooth projective surface with $\omega_X \cong \OO_X$, $b_2(X) = 10$ and $\Pic^\tau_X \cong \alpha_2$. The associated $\alpha_2$-torsor induces a global $1$-form on $X$ and hence $h^0(X,T_X) = h^0(X,\Omega_X) = 1$. The next example proves Corollary \ref{Enriquescorollary}.

\begin{Example}[\emph{The automorphism scheme of generic supersingular Enriques surfaces}]\label{multiplicativeEnriques}
Assume that $X$ is generic. Then, it is known that $X$ contains no $(-2)$-curves (see e.g. \cite[Proposition 5.2]{Martinunnodal}) and that $X$ admits a multiplicative $p$-closed global vector field (see \cite[Theorem 8.16]{EkedahlHylShep}). By \cite[Theorem 5.7.1.]{CossecDolgachev}, there is an elliptic fibration $f: X \to \PP^1$, which, by \cite[Theorem 5.7.2.]{CossecDolgachev}, admits a unique multiple fiber $2F$, which is either additive or supersingular. By Remark \ref{vectorfield}, the existence of a multiplicative vector field implies $\mu_2 \subseteq \Aut_X^0$. Since $X$ contains no $(-2)$-curves, the fibration $f$ admits no reducible fibers, so Theorem \ref{mainrational} shows that $\im(f_*)$ is trivial. Hence, Lemma \ref{connsub} implies that $\Aut_X^0 \cong \Ker(f_*) \cong \mu_{2^n}$ and finally Theorem \ref{Main} shows that $n = 1$.
\end{Example}

\begin{Example}[\emph{$\Aut_X^0$-$1$-movable fibers exist on non-isotrivial surfaces}]
Again, let $X$ be a generic supersingular Enriques surface. In particular, we have $\Aut_X^0 = \mu_2$. By \cite[Theorem 3.4.1.]{CossecDolgachev}, the surface $X$ admits two elliptic fibrations $f_1,f_2: X \to \PP^1$ with unique double fibers $2F_i$ satisfying $F_1.F_2 = 1$. We claim that $F_1$ is $\Aut_X^0$-movable if it is smooth. In fact, one can show that this condition is automatically satisfied for generic $X$, but for the sake of brevity we will not prove this here.
Seeking a contradiction, we assume that $F_1$ is not $\Aut_X^0$-$1$-movable. Since $F_1$ is smooth, it is supersingular, and thus it is fixed pointwise by $\Aut_X^0 = \mu_2$. But then $\mu_2$ fixes a point on a general fiber of $f_2$ and hence it fixes a general fiber of $f_2$ pointwise by Lemma \ref{fixedloci}. This is a contradiction and therefore $F_1$ is $\Aut_X^0$-$1$-movable.
\end{Example}

\begin{Remark}
Taking base changes of Example \ref{multiplicativeEnriques} along suitable finite and separable covers $C \to \PP^1$, one can construct non-isotrivial surfaces with non-trivial $\Ker(f_*)$ over every curve $C$ in characteristic $2$. We do not know how to construct similar examples if $p^n$ is bigger than $2$. This has essentially two reasons: First, the bounds given in Theorem (C) become very strong for $p^n \gg 0$ and second, it seems to be a very hard problem to construct elliptic surfaces with multiple supersingular and additive fibers whose multiplicity is a big power of $p$ (see for example \cite{Kawazoe} where this problem is studied in a very special case).
\end{Remark}

\subsection{Examples with non-trivial horizontal component} \label{sectionhorizontal}
In this section, we give examples of elliptic surfaces $f: X \to \PP^1$ where $\im(f_*)$ is non-trivial. More precisely, we will realize all cases described in Theorem \ref{mainrational} and the numbering of the examples will refer to the numbering in Theorem \ref{mainrational}. Since our examples admit a section, they will also show that all cases described in Theorem (D) occur.
 We will use the following technical Lemma, which allows us to construct some $\alpha_p$- or $\mu_p$-actions on $X$ by describing them on an affine Weierstrass equation.

\begin{Lemma} \label{extendvectorfield}
Let $f: X \to \PP^1$ be a Jacobian elliptic surface and let $D$ be a rational $p$-closed vector field on $X$.
Assume that $D$ is regular away from a fiber $F$ of $f$  and such that the induced rational vector field on $\PP^1$ is regular everywhere and has a zero at $f(F)$. Let $W$ be the isolated part of the zero locus of $D$, let $Z$ be the divisorial part, and let $S$ be a non-empty set of disjoint sections of $f$ to which $D$ is tangent. Then, $D$ is regular everywhere in each of the following cases:
\begin{enumerate}[(i)]
\item $F$ is of type ${\rm II}$ and ${\rm length}(W|_{X-F}) - (Z|_{X-F})^2 > c_2(X) - 4|S|$.
\item $p = 2$, $|S| \geq 2$, $F$ is of type ${\rm III}$, ${\rm length}(W|_{X-F}) - (Z|_{X-F})^2 > c_2(X) - 6$, and $D^2 = D$.
\end{enumerate}
\end{Lemma}

\prf
Let $t = f^{\#}(s)$, where $s$ is a parameter at $f(F)$. Let $S = \{\Sigma_1,\hdots,\Sigma_n\}$ be disjoint sections of $f$ such that $D$ is tangent to $\Sigma_i$ and let $x_i$ be a local equation for $\Sigma_i$ in a neighborhood of $P_i := F \cap \Sigma_i$. Then, in the completion $\widehat{\OO}_{X,P_i} \cong k[[x_i,t]]$, we can write $D$ as
$$
D = t^{-{m_i}}(t^{m_i+l} \frac{\partial}{\partial t} + g_i \frac{\partial}{\partial x_i}),
$$
where $m_i$ is the pole order of $D$ along the component of $F$ meeting $\Sigma_i$, $l \geq 1$ is the zero order of the induced vector field on $\PP^1$ at $f(F)$ and $g_i \in k[[x_i,t]]$ is a power series with $t \nmid g$ and $x_i \mid g$, since $D$ is tangent to $\Sigma_i$. 
In particular, $W$ has multiplicity ${\rm mult}_{P_i}(W) = (m_i+l) \cdot {\rm mult}_{P_i}(g_i)$ at $P_i$.

Since $Z$ is contained in fibers of $f$, we can apply \cite[Proposition 2.1]{KatsuraTakeda} to the part $W'$ of $W$ with support in $F$ to obtain
\begin{equation}\label{eqnvectorfield}
c_2(X) - \sum_{i=1}^n {\rm mult}_{P_i}(W) \geq c_2(X) - {\rm length}(W') =  {\rm length}(W|_{X-F}) - Z^2. \tag{$\ast$}
\end{equation}

Assume first that $F$ is of type ${\rm II}$ and $D$ has a pole along $F$. Then, all the $m_i$ are equal and $m := m_1 > 0$. Moreover, $t^mD$ is a regular $p$-closed vector field near $F$. We have $(t^mD)^p(t) = 0$, hence $t^mD$ is additive and thus ${\rm mult}_{P_i}(g_i) \geq 2$. But then ${\rm mult}_{P_i}(W) = (m+l) \cdot {\rm mult}_{P_i}(g_i) \geq 4$. Plugging into equation \eqref{eqnvectorfield}, this proves Claim (i).

Next, assume that $F$ is of type ${\rm III}$. If $D$ has poles along both components $F_1$ and $F_2$ of $F$, then Claim (ii) follows by the same argument as in the previous paragraph. If $D$ has a pole along $F_1$ but not along $F_2$, then we consider the contraction $\pi:X \to X'$ of $F_1$. Then, $D$ induces a $\mu_2$-action on $X'$ that preserves the image $F'$ of $F$. Note that $F'$ is a cuspidal rational curve. Moreover, by Proposition \ref{blowup}, the $\mu_2$-action does not fix the cusp of $F'$, for otherwise it would lift to $X$. Hence, by Lemma \ref{fixedloci}, the $\mu_2$-action on $F'$ has four isolated fixed points. In particular, we have ${\rm length}(W') \geq 4$. On the other hand, the pole of $D$ along $F_1$ contributes at least $(-2)$ to the right hand side of Equation \eqref{eqnvectorfield}. Hence,
$$
{\rm length}(W|_{X-F}) - (Z|_{X-F})^2 \leq {\rm length}(W|_{X-F}) - Z^2 - 2 \leq c_2(X) - {\rm length}(W') - 2 \leq c_2(X) - 6
$$
contradicting our assumption. This proves Claim (ii).
\qed
\vspace{3mm}

\begin{Example}[\emph{Case $(i)$}]\label{casea}
Consider the following four affine Weierstrass equations, where $u,v \in k$ are parameters and $t$ is a coordinate on $C = \PP^1$:
\begin{eqnarray*}
y^2 &=& x^3  + t \\
y^2 &=& x^3 + tx \\
y^2 &=& x^3 + t^2 \\ 
y^2 &=& x^3 + ut^2x + vt^3
\end{eqnarray*}
The induced minimal proper regular models $f: X \to \PP^1$ are precisely the four types of surfaces described in Lemma \ref{twofibers} (see \cite{MirandaPersson}). Now, note that each of these Weierstrass models admits a $\GG_m$-action given by $t \mapsto \lambda^a t, x \mapsto \lambda^2 x, y \mapsto \lambda^3 y$, where $a = 6,4,3,$ and $2$, respectively. Moreover, since $\GG_m$ is smooth and $X$ is the minimal resolution of the corresponding Weierstrass model, we obtain a $\GG_m$-action on $X$. Since $\Ker(f_*) \cap \GG_m$ is finite in every case, we have $\GG_m \subseteq \im(f_*)$.
\end{Example}

\begin{Remark}
Note that for the first surface, the $\GG_m$-action on $X$ induces the vector field $D = 6t\frac{\partial}{\partial_t} + 2x\frac{\partial}{\partial_x}$ in a neighborhood of the fiber of type ${\rm II}$ at $t = 0$. This $D$ is a counterexample to \cite[Lemma 4]{RudakovShafarevich} for all $p > 3$. The problem with the proof of \cite[Lemma 4]{RudakovShafarevich} is that not every vector field on a Weierstrass model is of the form claimed there.
\end{Remark}

\begin{Example}[\emph{Case $(ii)$ with $p = 3$}]\label{caseb3}
Let $p =3$ and consider the following three affine Weierstrass equations, where $k \geq 1$ is an integer and $t$ is a coordinate on $C = \PP^1$:
\begin{eqnarray*}
y^2 &=& x^3 + tx + t \\
y^2 &=& x^3 + x^2 + t^{3^{2k}}\\
y^2 &=& x^3 + tx^2 + t^{3^{2k-1} + 3} 
\end{eqnarray*}
We claim that that the corresponding elliptic surface $f: X \to \PP^1$ has precisely two fibers of type $({\rm II},{\rm III}), ({\rm II},{\rm I}_{3^{2k}}),$ and  $({\rm II},{\rm I}_{3^{2k-1}}^*)$, respectively and that $X$ admits a $\mu_3$-action which is non-trivial on $\PP^1$ in each of these cases.

The first Weierstrass model $X'$ can be embedded in $\PP(1,1,2,3)$ as
$$
y^2 = x^3 + ts^3 x + ts^5
$$
and it follows immediately from Tate's algorithm that $f$ admits a fiber of type ${\rm II}$ over $t = 0$ and a fiber of type ${\rm III}^*$ over $s = 0$. Note that the surface admits a unique singularity at the point $P$ given by $[s:t:x:y] = [0:1:0:0]$ and this singularity is a rational double point of type $E_7$. There is a $\mu_3$-action on the Weierstrass model given by 
$$
[s:t:x:y] \mapsto [s:at:a^2x + (1-a)s^2:y] \hspace{1cm} a^3 = 1.
$$
This action fixes $P$, so it lifts to the blow-up of $X'$ at $P$ by Proposition \ref{blowup}. By \cite[Theorem 4.1 (iii)]{Hirokado}, this already implies that the $\mu_3$-action lifts to $X$.

The second Weierstrass model $X'$ is an affine chart of the pullback along the $(2k)$-fold Frobenius on $\PP^1$ of the surface $Y' \subseteq \PP(1,1,2,3)$ given by
$$
y^2 = x^3 + s^2x^2 + s^5t.
$$
By Tate's algorithm, the minimal proper regular model $g: Y \to \PP^1$ of $Y'$ has a fiber $F'$ of type ${\rm II}^*$ with $\delta_{F'} = 1$ over $s = 0$ and a fiber of type ${\rm I}_1$ over $t = 0$. Since the Swan conductor does not change if we pull back along Frobenius and the vanishing order of $\Delta_g$ gets multiplied by $3$, the elliptic surface $f: X \to \PP^1$ admits a fiber $F$ over $s = 0$ with $\delta_{F} = 1$ and $v_{f(F)}(\Delta_f)= 11 \cdot 3^{2k} = 3$ mod $12$. This shows that $F$ is of type ${\rm II}$. Moreover, the fiber of $f$ over $t = 0$ is of type ${\rm I}_{3^{2k}}$. 
The affine Weierstrass equation for $X$ admits a $\mu_3$-action given by
$$
(t,x,y) \mapsto (at,x,y) \hspace{1cm} a^3 = 1.
$$
As in the previous case, this $\mu_3$-action preserves the singular point $(0,0,0)$ of the affine Weierstrass equation and lifts to the minimal resolution. Moreover, the $\mu_3$-action corresponds to the rational vector field $D = t\frac{\partial}{\partial_t}$ on $X$ which is regular away from $F$ and tangent to the zero section of $f$. A straightforward local computation shows that  ${\rm length}(\langle D \rangle|_{X-F}) - ((D)|_{X-F})^2 = 3^{2k} = c_2(X) - 3$, hence $D$ is regular on $X$ by Lemma \ref{extendvectorfield}, giving the desired $\mu_3$-action on $X$.

The third Weierstrass model $X'$ is an affine chart of the quadratic twist by $t$ of the pullback along the $(2k-1)$-fold Frobenius on $\PP^1$ of the surface $Y' \subseteq \PP(1,1,2,3)$ given by
$$
y^2 = x^3 + s^2x^2 + s^5t.
$$
By a similar argument as in the previous case, the minimal proper regular model of the pulled back surface admits a singular fiber of type ${\rm II}^*$ and a fiber of type ${\rm I}_{3^{2k-1}}$. Therefore, by Lemma \ref{twist} (i), the singular fibers of the minimal proper regular model $f: X \to \PP^1$ of $X'$ are of the stated types. There is a $\mu_3$-action on the affine chart $X'$ given by
$$
(t,x,y) \mapsto (at,ax,y) \hspace{1cm} a^3 = 1.
$$
The rest of the argument is similar to the previous case.
\end{Example}

\begin{Remark}
We remark that all three of the above surfaces are counterexamples to \cite[Lemma 4]{RudakovShafarevich} in characteristic $3$. Moreover, the second and third example are counterexamples to \cite[Theorem 6]{RudakovShafarevich}. The proof of this Theorem fails in Case (6), where \cite[Lemma 4]{RudakovShafarevich} is applied. Moreover, note that our equation for the fibration with fibers of type $({\rm II},{\rm I}_{3^{2k-1}}^*)$ differs from Equation (3) given in \cite{RudaShafaVector} and the corresponding equation given in \cite[p.1503]{RudakovShafarevich2}. Using Tate's algorithm, one can check that, at least for general $k$, these two equations do not admit a fiber $F$ of type ${\rm II}$ with $v_{f(F)}(\Delta_f) = 3$, hence they do not admit global vector fields.
\end{Remark}

\begin{Example}[\emph{Case $(ii)$ with $p = 2$}]\label{caseb2} Let $p =2$ and consider the following four affine Weierstrass equations, where $k \geq 1$ is an integer and $t$ is a coordinate on $C = \PP^1$:
\begin{eqnarray*}
y^2 + ty &=& x^3 + t \\
y^2 + ty &=& x^3 \\
y^2 + xy &=& x^3 + t^{2^{2k}}x \\
y^2 + xy &=& x^3 + t^{2^{2k-1}}x^2 + t^{2^{2k}}x
\end{eqnarray*}
We claim that that the corresponding elliptic surface $f: X \to \PP^1$ has precisely two fibers of type $({\rm II},{\rm IV}^*),({\rm III},{\rm IV}^*), ({\rm III},{\rm I}_{2^{2k+1}}),$ and  $({\rm II},{\rm I}_{2^{2k+1}})$, respectively, and that $X$ admits a $\mu_2$-action which is non-trivial on $\PP^1$ in every case.

The first Weierstrass model $X'$ can be embedded in $\PP(1,1,2,3)$ as 
$$
y^2 + s^2ty = x^3 + s^5t
$$
and it follows from Tate's algorithm that $f$ admits a fiber of type ${\rm II}$ over $t = 0$ and a fiber of type ${\rm IV}^*$ over $s = 0$. There is a $\mu_2$-action on the Weierstrass model given by
$$
[s:t:x:y] \mapsto [as:t:x:y + (1+a)s^3] \hspace{1cm} a^2 = 1.
$$
This action fixes the unique singular point $P = [0:1:0:0]$ of $X'$, hence it lifts to the blow-up of $X'$ at $P$ by Proposition \ref{blowup}. Since $P$ is of type $E_6$, it follows from \cite[Theorem 5.1 (iii)]{Hirokado} that the action lifts to $X$.

Similarly, the second Weierstrass model $X'$ can be embedded in $\PP(1,1,2,3)$ as 
$$
y^2 + s^2ty = x^3.
$$
This time, Tate's algorithm shows that $f$ admits a fiber of type ${\rm III}$ over $t = 0$ and a fiber of type ${\rm IV}^*$ over $s = 0$. There is a $\mu_2$-action on the Weierstrass model given by 
$$
[s:t:x:y] \mapsto [as:t:x:y] \hspace{1cm} a^2 = 1.
$$
This action fixes the two singular points $[0:1:0:0]$ and $[1:0:0:0]$ and hence, as in the previous case, it lifts to $X$.

The third Weierstrass model $X'$ is an affine chart of the pullback along the $(2k-2)$-fold Frobenius on $\PP^1$ of the surface $Y' \subseteq \PP(1,1,2,3)$ given by
$$
y^2 + sxy = x^3 + t^4x.
$$
By Tate's algorithm, the minimal proper regular model $g: Y \to \PP^1$ of $Y'$ has a fiber $F'$ of type ${\rm III}$ with $\delta_{F'} = 1$ over $s = 0$ and a fiber of type ${\rm I}_8$ over $t = 0$. Similarly to the analogous case if $p = 3$, it is easy to check that $f:X \to \PP^1$ admits a fiber $F$ of type ${\rm III}$ with $\delta_F = 1$ over $s = 0$ and a fiber of type ${\rm I}_{2^{2k+1}}$ over $t = 0$. There is a $\mu_2$-action on $X'$ given by
$$
(t,x,y) \mapsto (at,x,y) \hspace{1cm} a^2 = 1.
$$
This action corresponds to the vector field $D = t \frac{\partial}{\partial t}$ on $X$, which satisfies $D^2 = D$ and is tangent to the zero section $\Sigma_1$ and to the $2$-torsion section $\Sigma_2$ given by $x = y = 0$. Using the height pairing (see \cite[p.110]{schuettShioda}), one can check that $\Sigma_2$ is disjoint from $\Sigma_1$. Since $D$ fixes $(0,0,0)$, it lifts to the minimal resolution of this singularity. Moreover, a straightforward local computation shows that ${\rm length}(\langle D \rangle|_{X-F}) - ((D)|_{X-F})^2 = 2^{2k+1} = c_2(X) - 4$, so that $D$ is regular on all of $X$ by Lemma \ref{extendvectorfield}. This yields the desired $\mu_2$-action on $X$.
 
The fourth Weierstrass model $X'$ is an affine chart of the pullback along the $(2k-2)$-fold Frobenius on $\PP^1$ of the surface $Y' \subseteq \PP(1,1,2,3)$ given by
$$
y^2 + sxy = x^3 + t^2x^2 + t^4x.
$$
By Tate's algorithm, the minimal proper regular model $g: Y \to \PP^1$ of $Y'$ has a fiber $F'$ of type ${\rm II}$ with $\delta_{F'} = 2$ over $s = 0$ and a fiber of type ${\rm I}_8$ over $t = 0$. Using Tate's algorithm, one can check that a fiber of type ${\rm II}$ with $\delta_{F'} = 2$ remains of the same type when pulled back along an even power of the Frobenius, hence $f: X \to \PP^1$ admits a fiber of type ${\rm II}$ with $\delta_{F} = 2$ over $s = 0$ and a fiber of type ${\rm I}_{2^{2k+1}}$ over $t = 0$. The $\mu_2$-action on $X'$ given by 
$$
(t,x,y) \mapsto (at,x,y) \hspace{1cm} a^2 = 1
$$
extends to a $\mu_2$-action on $X$ by the same argument as in the previous case.
\end{Example}

\begin{Remark}
We remark that the first and the fourth of the above surfaces are counterexamples to \cite[Lemma 4]{RudakovShafarevich} in characteristic $2$ and the second and the third are counterexamples to \cite[Lemma 3]{RudakovShafarevich}. Moreover, if we choose $k$ such that $2^{2k+1}+4 = 12$ mod $24$, we obtain counterexamples to \cite[Theorem 6]{RudakovShafarevich}. Again, the proof of the latter fails in Case (6), where the erroneous Lemmas 3 and 4 are applied. Moreover, we remark that the surfaces with fibers of type $({\rm III},{\rm I}_{2^{2k+1}})$ are missing from the classification in \cite{RudaShafaVector}.
\end{Remark}

\begin{Example}[\emph{Cases $(iii)$ and $(iv)$}]\label{casecd}
Consider the following affine Weierstrass equations, where $t$ is a coordinate on $C = \PP^1$ and $u \in k^\ast$:
\begin{eqnarray*}
p=3: & y^2 &= x^3 + x + t \\
p=2: & y^2 + y &= x^3 + t \\
& y^2 + uxy &= x^3 + tx^2 + x
\end{eqnarray*}
Each of these surfaces admits a $\GG_a$-action given by
\begin{eqnarray*}
(t,x,y) &\mapsto& (t + a^3 + a, x - a,y) \\
(t,x,y) &\mapsto& (t + a^2 + a, x,y + a) \\
(t,x,y) &\mapsto& (t + a^2 + ua, x, y + ax) \hspace{1cm} a \in k
\end{eqnarray*}
which lifts to the respective minimal proper regular model $f: X \to \PP^1$. The first two surfaces admit a unique singular fiber of type ${\rm II}^*$ over $t = \infty$ and the generic fiber of $f$ is supersingular. The third surface admits a unique singular fiber of type ${\rm I}_4^*$  over $t = \infty$ and the generic fiber of $f$ is ordinary with $j$-invariant $u^8$.
\end{Example}

\section{Proofs of the Main Theorems}\label{proofs}
In this section, we combine our study of horizontal and vertical components of $\Aut_X^0$ in order to prove Theorem (A), (B), (C), and (D) of the introduction. Moreover, we recall how the non-existence of global vector fields on K3 surfaces follows from Theorem (D).
\vspace{3mm}

{\sc Proof of Theorem (A)}
\vspace{1mm}

Let us prove Claim (i). Since $f$ is not isotrivial, we have $\Ker(f_*) \cong \mu_{p^n}$ for some $n \geq 0$ by Lemma \ref{connsub}. Moreover, by Proposition \ref{Propelliptic}, Proposition \ref{Proprational} and Theorem \ref{mainrational}, the horizontal component $\im(f_*)$ is trivial unless possibly in the cases described in Theorem \ref{mainrational} (ii). In these latter cases, we have $\Aut_X^0 \cong \im(f_*) \subseteq \mu_p$, so $h^0(X,T_X) \leq 1$ holds in every case.

As for Claim (ii), assume that the generic fiber of $f$ is ordinary or that $f$ admits no multiple fibers, and that $h^0(X,T_X) \geq 2$. Then, by Lemma \ref{connsub} and Corollary \ref{corollarymain}, we have $\Ker(f_*)^0 \in \{\mu_{p^n},M_n,E\}$ where $n \geq 0$ and $E$ is an elliptic curve, so $\im(f_*)$ has to be non-trivial. Now, Proposition \ref{Propelliptic}, Proposition \ref{Proprational} and Theorem \ref{mainrational} imply that $X$ is either ruled over an elliptic curve, an Abelian surface isogeneous to a product of elliptic curves, bielliptic with $\omega_X \cong \OO_X$, or an elliptic surface $f: X \to \PP^1$ with a unique singular fiber, without multiple fibers, and with supersingular generic fiber. In the first case, we have described the automorphism scheme in Example \ref{ruled}. In particular, we have seen that $h^0(X,T_X) \leq 4$ holds. In the second and third case, we have $h^0(X,T_X) = 2$ by Proposition \ref{Propelliptic}. In the fourth case, the vertical component $\Ker(f_*)^0$ is trivial by Theorem \ref{Main} and the horizontal component $\im(f_*)$ fixes a point on $\PP^1$ and thus $\im(f_*) \subseteq \GG_a \rtimes \GG_m$. In particular, we have $\im(f_*)[F] \subseteq \GG_a \rtimes \GG_m[F] = \alpha_p \rtimes \mu_p$. Now, Theorem \ref{mainrational} shows that $\im(f_*)[F] = \alpha_p$, hence $h^0(X,T_X) \leq 1$, so this case does not occur.

Finally, Claim (iii) is Example \ref{importantexample}.
\qed
\vspace{3mm}

{\sc Proof of Theorem (B)}
\vspace{1mm}

Assume first that $\im(f_*)$ is non-trivial. Then, by Proposition \ref{Propelliptic}, Proposition \ref{Proprational} and Theorem \ref{mainrational}, we have $\Aut_X^0 \cong \im(f_*) \cong \mu_p$ with $p \in \{2,3\}$. If $\im(f_*)$ is trivial, then $\Aut_X^0 \cong \Ker(f_*) \cong \mu_{p^n}$ for some $n \geq 0$ by Lemma \ref{connsub}. This proves Theorem (B).
\qed
\vspace{3mm}

{\sc Proof of Theorem (C)}
\vspace{1mm}

The inequality is trivial if $p^n \in \{2,3\}$, so we may assume $p^n \geq 4$. Then, by Proposition \ref{Propelliptic}, Proposition \ref{Proprational} and Theorem \ref{mainrational}, we have $\Aut_X^0 \cong \Ker(f_*) \cong \mu_{p^n}$ and then the inequality is proved in Proposition \ref{Igusainequality}. Next, note that by the results of Section \ref{horizontal}, the conditions given in Theorem (C) guarantee that $\Aut_X^0 \cong \Ker(f_*)$, so the statement about the multiplicities of additive and supersingular fibers of $f$ is exactly Theorem \ref{Main} (iii) (2).
\qed
\vspace{3mm}

{\sc Proof of Theorem (D)}
\vspace{1mm}

Assume that $c_2(X) \neq 0$. Since $f$ admits no multiple fibers, Theorem \ref{Main} shows that $\ker(f_*)^0$ is trivial. Moreover, by Proposition \ref{Propelliptic} and Proposition \ref{Proprational}, we have $C = \PP^1$ and $b_1(X) = 0$. In particular, $f$ admits a singular fiber and thus we have an inclusion $\Aut_X^0 \cong \im(f_*) \subseteq \GG_a \rtimes \GG_m$, since $\GG_a \rtimes \GG_m$ is the stabilizer of a point on $\PP^1$. Thus, either $\alpha_p \subseteq \im(f_*)$ or $\mu_p \subseteq \im(f_*)$.

If $\alpha_p \subseteq \im(f_*)$, then the $\alpha_p$-action on $X$ preserves every singular fiber of $f$ by Remark \ref{criteriaremark} and only one point on $\PP^1$ by Lemma \ref{fixedloci}, hence $f$ is isotrivial with a unique singular fiber.  In particular, by Theorem \ref{mainrational} (iv), we have $p \in \{2,3\}$. Thus, we are in Case (v) of Theorem (D).

If $\mu_p \subseteq \im(f_*)$, then Theorem \ref{mainrational} shows that either the singular fibers of $f$ are as in Lemma \ref{twofibers} or $p \in \{2,3\}$ and the singular fibers of $f$ are as in Lemma \ref{twofibers23}. In particular, if the fibers are not of the types described in Theorem (D) (iii) and (iv), then $f$ is isotrivial and $X$ satisfies $c_2(X) = 12$
by Ogg's formula. Hence $\chi(X,\OO_X) = 1$ and therefore $\omega_X \cong \OO_X(-F)$, where $F$ is the class of a fiber of $f$. In particular, we have $h^1(X,\OO_X) = h^2(X, \omega_X^{\otimes 2}) = 0$ and thus $X$ is rational and $f$ admits a section. This is Case (ii) of Theorem (D).
\qed
\vspace{3mm}

{\sc Proof of Corollary \ref{K3corollary}}
\vspace{1mm}

Let $X$ be a K3 surface and assume by contradiction that $h^0(X,T_X) \neq 0$. By \cite[p. 1502]{RudakovShafarevich2}, this implies that the surface $X$ admits an elliptic fibration $f: X \to \PP^1$ with at least two singular fibers. This contradicts Theorem (D), because no elliptic surface listed in Theorem (D) (i)-(iv) satisfies the equality $c_2(X) = 24$, which holds for the K3 surface $X$. Therefore, we must have $h^0(X,T_X) = 0$.
\qed

\newpage

\section{Appendix: Some quadratic twists} \label{appendix}
In this section, we give some background on quadratic twists, which we needed for example in the proof of Lemma \ref{twofibers}. Let $f: X \to C$ be a Jacobian elliptic surface. Then, a \emph{quadratic twist} of $f$ is a Jacobian elliptic surface $f': X' \to C$  that becomes isomorphic to $f$ after passing to a degree two cover of $C$. If the generic fiber of $f$ is ordinary, then all its twists are quadratic. To make this more explicit, let $d \in k(C)$ be a rational function. Then, the quadratic twist $f_d: X_d \to C$ of $f$ by $d$ is defined as follows: If $p \neq 2$ and the generic fiber of $f$ is given by
$$
y^2 = x^3 + a_2x^2 + a_4x + a_6
$$
with $a_i \in k$, then $f_d$ is given by
$$
y^2 = x^3 + da_2x^2 + d^2a_4x + d^3a_6.
$$
If $p = 2$ and the generic fiber of $f$ is given by
$$
y^2 + a_1xy + a_3y = x^3 + a_2x^2 + a_4x + a_6
$$
with $a_i \in k$, then $f'$ is given by 
$$
y^2 + a_1xy + a_3y = x^3 + (a_2 + da_1^2)x^2 + a_4x + a_6 + da_3^2.
$$
The fibers of $f$ and $f_d$ are isomorphic except possibly over the set $S$ of poles and zeroes of $d$ (resp. the set of poles if $p = 2$) and we say that $f_d$ is a \emph{quadratic twist of $f$ at $S$}. Quadratic twists by $d_1$ and $d_2$ are isomorphic if and only if $d_1/d_2$ is a square if $p \neq 2$ (resp. if and only if $d_1 + d_2 = c^2 + c$ for some $c \in k(C)$ if $p = 2$).
In the following lemma, we summarize the facts about quadratic twists that we used in this article.
\begin{Lemma} \label{twist}
Let $f: X \to C$ be an elliptic surface and let $d \in k(C)$. Let $F$ be a fiber of $f$ and $F_d$ the corresponding fiber of $f_d$. Then, the following hold:
\begin{enumerate}[(i)]
\item If $p \neq 2$ and $d$ has a zero or pole at $f(F)$, then the types of $F$ and $F_d$ are related as follows: ${\rm I}_n \leftrightarrow {\rm I}_n^*$, ${\rm II} \leftrightarrow {\rm IV}^*, {\rm III} \leftrightarrow {\rm III}^*, {\rm IV} \leftrightarrow {\rm II}^*$.
\item If $p = 2$, $C = \PP^1$ and $F$ is of type ${\rm II}$ with $\delta_{F} = 2$, then we can choose $d \in k(t)$ with a single simple pole such that $F_d$ is of type ${\rm III}$ with $\delta_{F} = 1$.
\item If $p = 2$ , $C = \PP^1$ and $f$ has ordinary generic fiber and a unique singular fiber $F$, then $f$ is isotrivial and $F$ is of type ${\rm I}_{8k + 4}^*$ with $\delta_{F} = 4k + 2$ for some $k \geq 0$.
\end{enumerate}
\end{Lemma}

\prf
Claim (i) is well-known, see for example \cite[Section 5.4.]{schuettShioda}.

As for Claim (ii), choose a Weierstrass equation
$$
y^2 + a_1xy + a_3y = x^3 + a_2 x^2 + a_4 x + a_6
$$
with coefficients $a_i \in k[t]$, where $t$ is a parameter at $f(F)$. Then, $F$ being of type ${\rm II}$ with $\delta_F = 2$ means that we can choose the $a_i$ such that $t \mid a_1,a_3,a_4,a_6$ but $t^2 \nmid a_3,a_6$. Let $c_3$ resp. $c_6$ be the linear terms of $a_3$ resp. $a_6$. If we set $d = c_6/c_3^2$, the quadratic twist
$$
y^2 + a_1xy + a_3y = x^3 + (a_2 + da_1^2)x^2 + a_4 x + a_6 + da_3^2
$$
still has coefficients in $k[t]$ and we have $t^2 \mid a_6 + da_3^2$. Note that $t^3 \nmid b_8 := (a_1^2a_6 + a_1a_3a_4 + a_2a_3^2 + a_4^2)$ and the quadratic twist does not change $b_8$. Thus, Tate's algorithm shows that $F_d$ is of type ${\rm III}$ with $\delta_{F_d} = 1$. Moreover, the twist parameter $d$ has a simple pole at $f(F)$ and no other poles.

Next, let us prove Claim (iii). First, we prove that $f$ is isotrivial. For this, choose a parameter $t$ on $\PP^1$ such that $F$ is located at $t = 0$. By \cite[Appendix A]{Silverman}, the assumption that the generic fiber of $f$ is ordinary allows us to find a Weierstrass equation of the form
$$
y^2 + xy = x^3 + \frac{a}{b} x^2 + a_6
$$
with $a,b \in k[t]$ and $a_6 \in k(t)$. Write $a/b = \sum_{i = -n}^\infty d_it^i \in k((t))$ and twist the above equation by $d = \sum_{i = -n}^{-1} d_it_i \in k(t)$. This quadratic twist only changes the fiber over $t = 0$, so we may assume that $t \nmid b$.

We have $\Delta = a_6$ and $j = 1/a_6$. Since $f$ has no singular fibers away from $t = 0$, the $j$-map has no poles away from $t = 0$ and $\Delta$ is constant up to $12$-th powers. Therefore, we can write $a_6 = t^{12n}/c^{12}$ for some $n \geq 0$ and $c \in k[t]$ with $t \nmid c$ and $\deg(c) \leq n$. Then, we can rescale the Weierstrass equation to an integral Weierstrass equation of the following form
$$
y^2 + bc^{2}xy = x^3 + abc^{4} x^2 + t^{12n}b^6.
$$
If $n > 0$, then Tate's algorithm shows that $F$ is of type ${\rm I}_{12n}$, because $t \nmid b,c$. Then, Igusa's inequality shows $12n \leq b_2(X) = c_2(X) - 2$, which contradicts Ogg's formula $12n = c_2(X)$. Hence, we must have $n = 0$ and thus $j$ is constant.

This implies that the generic fiber of $f$ is a quadratic twist of the ordinary elliptic curve with $j$-invariant $j$ given by
$$
y^2 + xy = x^3 + j
$$
by a twist parameter $d \in k(t)$ whose only poles are at $t = 0$.
Every non-trivial such twist can be written as
$$
y^2 + xy = x^3 + \frac{a}{t^{2k+1}}x^2 + j
$$
for some $a \in k[t]$ of degree at most $2k+1$ with $t \nmid a$, where $k \geq 0$ is an integer. Clearing denominators and applying $y \mapsto \sqrt{j} t^{3k + 3}$, we obtain the equation
$$
y^2 + t^{k+1}xy = x^3 + atx^2 + \sqrt{j}t^{4k+4}x.
$$
By Tate's algorithm, this equation is minimal and $F$ is of type ${\rm I}_{8k+4}^*$. Moreover, $\Delta = t^{12k+12}$, so that $\delta_{F} = 4k + 2$.
\qed

\bibliographystyle{alpha} 
\bibliography{AutomorphismsEllipticSurfaces}
\end{document}